\documentclass[a4paper,10pt]{article}
\usepackage[margin=1.5in]{geometry}
\usepackage{lipsum}
\usepackage{graphicx}
\usepackage{epstopdf}
\usepackage{algorithmic}
\ifpdf
  \DeclareGraphicsExtensions{.eps,.pdf,.png,.jpg}
\else
  \DeclareGraphicsExtensions{.eps}
\fi

\usepackage{graphicx}
\usepackage{hyphenat}
\usepackage[dvipsnames]{xcolor}
\usepackage{amssymb, amsmath, amsfonts, array}
\usepackage{mathtools, xfrac, scalerel}
\usepackage[export]{adjustbox}

\usepackage[linesnumbered, ruled, noresetcount, vlined]{algorithm2e}
\SetKwFor{Parfor}{parfor}{do}{}

\usepackage{tikz}
\usepackage{pgfplots}
\usepackage{pgfplotstable}
\usepgfplotslibrary{patchplots}
\pgfplotsset{compat=1.15}
\usepackage{arydshln}
\usetikzlibrary{math, matrix, shapes, calc, arrows.meta, decorations.pathreplacing, decorations.markings}
\usetikzlibrary{colorbrewer}
\usepackage{listofitems,siunitx,enumerate}

\usepackage[hyphens]{url}
\usepackage[hidelinks]{hyperref}
\usepackage[nocompress]{cite}

\usepackage{etoolbox}
\apptocmd{\sloppy}{\hbadness 10000\relax}{}{}

\newcommand{\secref}[1]{Sec.~\ref{#1}}
\newcommand{\figref}[1]{Fig.~\ref{#1}}
\newcommand{\algref}[1]{Alg.~\ref{#1}}
\newcommand{\verteqiv}{\rotatebox{90}{$\,\equiv$}}
\newcommand{\equivto}[2]{\underset{\scriptstyle\overset{\mkern4mu\verteqiv}{#2}}{#1}}
\DeclarePairedDelimiter\floor{\lfloor}{\rfloor}

\newif\ifcompileFigs

\title{Multigrid Reduction in Time for non-linear hyperbolic equations\date{}\thanks{This publication is based on work partially supported by the EPSRC Centre For Doctoral Training in Industrially Focused Mathematical Modelling (EP/L015803/1) in collaboration with the Culham Centre for Fusion Energy, and by the National Productivity Investment Fund (NPIF). The work of S. MacLachlan was partially supported by NSERC discovery grants RGPIN-2014-06032 and RGPIN-2019-05692.}
}

\author{Federico Danieli\thanks{Mathematical Institute, University of Oxford, Oxford, UK 
  (\href{mailto:federico.danieli@maths.ox.ac.uk}{\tt{federico.danieli@maths.ox.ac.uk}}.)},
   \and Scott MacLachlan\thanks{Department of Mathematics and Statistics, Memorial University of Newfoundland, St. John's, NL, Canada
  (\href{mailto:smaclachlan@mun.ca}{\tt{smaclachlan@mun.ca}}).}
}

\newenvironment{@abs}[1]{%
       \vspace{4pt}\footnotesize  \parindent 15pt {\bfseries #1. }\ignorespaces
     }
     {\par\vspace{7pt}}

\renewenvironment{abstract}{\begin{@abs}{\abstractname}}{\end{@abs}}
\newenvironment{keywords}{\begin{@abs}{\keywordsname}}{\end{@abs}}
\newenvironment{AMS}{\begin{@abs}{\AMSname}}{\end{@abs}}

\newcommand\keywordsname{Key words}
\newcommand\AMSname{AMS subject classifications}

\begin{document}

\maketitle

\begin{abstract}

Time-parallel algorithms seek greater concurrency by decomposing the temporal domain of a \emph{Partial Differential Equation} (PDE), providing possibilities for accelerating the computation of its solution.
While parallelisation in time has allowed remarkable speed-ups in applications involving parabolic equations, its effectiveness in the hyperbolic framework remains debatable: growth of instabilities and slow convergence are both strong issues in this case, which often negate most of the advantages from time-parallelisation.
Here, we focus on the \emph{Multigrid Reduction in Time} (MGRIT) algorithm, investigating in detail its performance when applied to non-linear conservation laws with a variety of discretisation schemes. Specific attention is given to high-accuracy \emph{Weighted Essentially Non-Oscillatory} (WENO) reconstructions, coupled with \emph{Strong Stability Preserving} (SSP) integrators, which are often the discretisations of choice for such PDEs. A technique to improve the performance of MGRIT when applied to a low-order, more dissipative scheme is also outlined.
This study aims at identifying the main causes for degradation in the convergence speed of the algorithm, and finds the \emph{Courant-Friedrichs-Lewy} (CFL) limit to be the principal determining factor.
\end{abstract}


\begin{AMS}{
  65M08,  	
  35L65,  	
  65M55,  	
  65Y05,  	
  65Y20   	
}
\end{AMS}

\section{Introduction}
\label{sec::intro}
\begin{figure}[!b]
  \centering
  \includegraphics[trim={0 0.52cm 0 1cm},clip,width=0.475\textwidth]{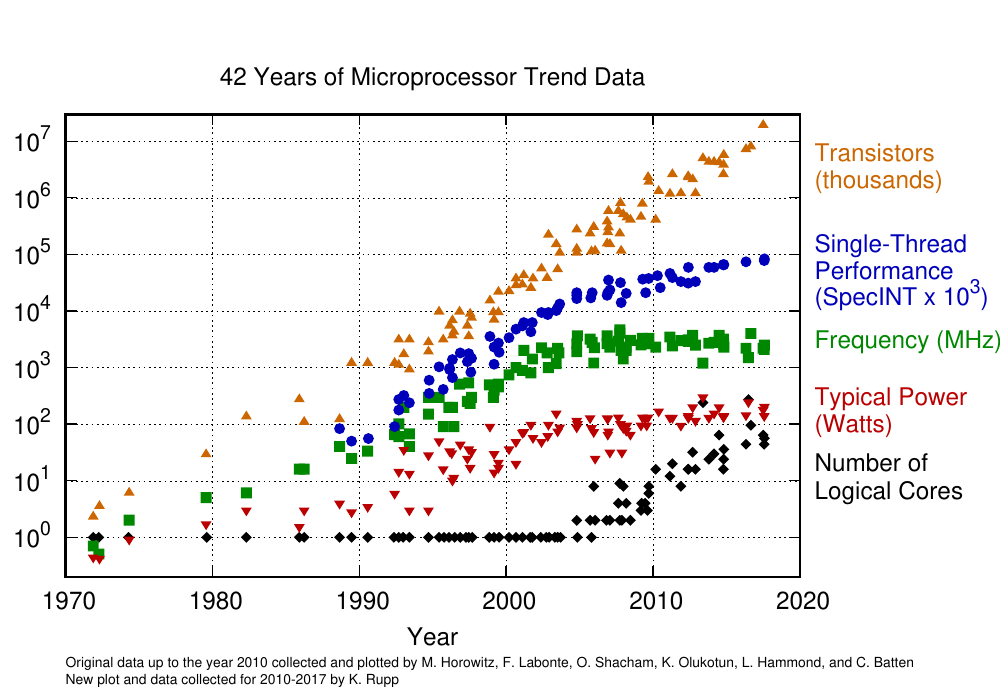}
  \caption[Evolution of microprocessors performance]{Evolution of microprocessors performance. In recent years, while the clock frequency of a single processor (green squares) has plateaued, the number of cores per processors (black diamonds) has been steadily increasing, which reflects the focus on parallelisation over power. Source: \cite[Creative Commons Attribution 4.0 International Public License]{processorPower}}
  \label{fig::introduction::processorPower}
\end{figure}

The growing complexity of computations arising from the numerical approximation of solutions to \emph{Partial Differential Equations} (PDEs) demands for ever-increasing computational power in order to tackle such problems in a reasonable amount of time. The frequency of computations that can be performed on a single processor unit, however, is limited (see \figref{fig::introduction::processorPower}), and represents an upper bound for the efficiency of serial algorithms. To overcome this limit, much attention has been directed towards the implementation of parallel algorithms, thus providing an alternative to sheer processing power to speed up the computation of such numerical solutions.
 
Some of the most standard parallelisation methods involve the decomposition of the spatial domain of the target PDE into smaller sub-domains, where some computations can be carried out independently \cite{ToselliWidlundDDM}. Information about solutions near interfaces needs to be exchanged between these sub-domains, usually in an iterative fashion. However, this procedure causes an additional overhead cost with respect to the serial solution, so that increasing the number of sub-domains above a certain limit becomes detrimental to the purpose of speeding up the computation, and spatial parallelisation is said to \emph{saturate}. Moreover, systems of \emph{Ordinary Differential Equations} (ODEs) do not present a spatial domain to subdivide in the first place.
For time-dependent PDEs and ODEs, the time domain presents an additional direction along which to seek parallelisation: if efficiently employed, this would allow us to make the most out of the parallel computation capabilities of modern supercomputers.

To achieve time-parallelisation, a variety of schemes have been introduced over the last few decades. Pioneering in this regard was the work by J-L. Lions, Y. Maday, and G. Turinici, who proposed the \emph{Parareal} algorithm in \cite{parareal}, but other methods quickly followed: among these, we point out the \emph{Parallel Full Approximation Scheme in Space and Time} (PFASST, \cite{PFASST,BoltenEtAl2018}), and the \emph{Revisionist Integral Deferred Correction} method (or RIDC, \cite{ChristliebEtAl2010,RIDC,OngEtAl2016}), as well as the \emph{Multigrid Reduction In Time} algorithm (MGRIT, \cite{MGRIToriginal,MGRIT_BDF,MGRITNonLin}), which is the main focus of this paper.  For a more thorough review of parallel-in-time methods, we refer to \cite{50yr}. Although different in their behaviour, at the core these schemes share a similar idea to achieve parallelisation: that is, pairing a \emph{fine integrator} (which is expensive to use, but is applied in parallel) together with a \emph{coarse integrator} (whose action is cheaper to compute, but is applied serially). The former is responsible for solving the target equation to the desired level of accuracy, while the latter takes care of quickly propagating updates along the time domain.
Time-parallelisation has proven very effective if used to speed up the solution of parabolic equations (see, for example, \cite{MGRITNonLin,pararealSkinTransport,pararealOptionPricing} for an analysis of some such applications from a range of different fields); unfortunately, though, its applicability to advection-dominated or hyperbolic PDEs still remains a matter of research.

Recovering an accurate numerical solution to hyperbolic equations \emph{per se} is challenging in many ways \cite[Chap.~1.4]{leveque}. One of the biggest difficulties lies in properly capturing shocks present in the solution, while retaining the desired level of accuracy and ensuring stability of the scheme used. These objectives are often at conflict with each other: better accuracy can be achieved using high-order approximations, but these approximations might generate spurious oscillations in proximity of discontinuities in the solution, undermining the stability of the algorithm. To counteract this behaviour, a variety of schemes have been proposed in the literature, some of which are described in this paper.
This problem is even more pronounced if time-parallelisation comes into play, when solutions from different integrators need to be combined together. Even if both fine and coarse solver are individually stable, often their combined use triggers instabilities that, in turn, result in loss of accuracy and poor convergence, negating most of the advantages from parallelisation. These issues were reported already in \cite{pararealAnalysis}, and emerge clearly from the analyses in \cite{pararealAnalysis2} and \cite{pararealAnalysis3}. Regardless, there remain some examples of successful applications: for example, in \cite{pararealHyp}, stability of Parareal applied to Burgers' equation is guaranteed by projecting the solution back to an energy manifold after every iteration; in \cite{pararealHyp2} the authors succeed in speeding up the solution to a hydrodynamic problem characterised by a Reynolds number of $5\cdot10^4$; in \cite{Nielsen}, choosing a Roe-average Riemann solver seems to be key for fast convergence when solving the shallow-water equations; in \cite{MGRIT_Oliver}, an optimisation approach is used to determine coarse-grid operators that achieve excellent performance of MGRIT for the linear advection equation.

In our work, we investigate in detail the performance of MGRIT when applied to some non-linear hyperbolic problems frequently employed as test-cases. In particular, we discuss the behaviour of the algorithm when used in conjunction with high-order discretisations, which are frequently used for the solution of conservation laws.
Example applications of parallel-in-time methods to high-order schemes are, to the knowledge of the authors, somewhat rare in the literature. Noticeable exceptions are: publications involving the PFASST algorithm \cite{PFASST}, which achieves high-order accuracy via repeated SDC corrections, but is limited to this very specific class of integrators; the works conducted in \cite{MGRIT_BDF,MGRIT_BDF2}, where $k$-th order multi-step BDF methods are analysed; \cite{Nielsen}, which shares a similar setup to ours, but only shows results for a specific configuration of the integrators used; finally, \cite{MGRIT_Oliver}, which is possibly the work whose goal is the closest to ours. There, the efficacy of MGRIT applied to hyperbolic equations is studied in detail, with the aim of identifying suitable coarse integrators to achieve fast convergence, though the project only covers the linear case.
The purpose of our study lies in investigating the applicability of MGRIT to the acceleration of solutions to non-linear systems of conservation laws which are relevant in real-world scenarios. It focuses on determining the principal factors that cause degradation in the convergence speed of the algorithm, and identifies future directions for improvements.

The remainder of this article is structured as follows.  In \secref{sec::modelProblems}, we present the test problems considered here and the discretisations used in our experiments. In \secref{sec::MGRIT}, we describe the MGRIT algorithm used in our simulations, and discuss the important choice of the coarse integrator for a given fine integrator. In \secref{sec::results}, we present and discuss numerical experiments.

\section{Model problems}
\label{sec::modelProblems}
As test cases for our experiments, we pick examples of \emph{conservation laws} that, due to their physical relevance, are widely used in the analysis and testing of numerical methods for hyperbolic equations. For simplicity, we limit ourselves to problems defined on a one-dimensional spatial domain, and with periodic boundary conditions. These equations can be written in the following general form:
\begin{equation}
  \left\{
  \begin{array}{lcr}
    \displaystyle\frac{\partial \boldsymbol{u}}{\partial t} + \frac{\partial \boldsymbol{f}(\boldsymbol{u})}{\partial x} = 0 &\quad\quad & (x,t)\in\left[0,L\right]\times\left[0,T\right]\\[\bigskipamount]
    \boldsymbol{u}(x,0) = \bar{\boldsymbol{u}}^0(x)   &  & x\in\left[0,L\right]\\
    \boldsymbol{u}(0,t) = \boldsymbol{u}(L,t)   &  & t\in\left[0,T\right]
  \end{array}
  \right.,
  \label{eqn::consLaw}
\end{equation}
where the variable $\boldsymbol{u}(x,t)\in\mathbb{R}^D$ contains the vector of conserved variables in the system, or \emph{state}, with associated \emph{flux} $\boldsymbol{f}(\boldsymbol{u})$. The parameters $L$ and $T$ define the size of the spatial and temporal domains, respectively, while $\bar{\boldsymbol{u}}^0(x)$ identifies the inital condition. We consider three different systems, ubiquitous in the literature.
\begin{itemize}
  \item The \emph{Burgers' equation}, here considered in its inviscid formulation, which is possibly the simplest example of a scalar conservation law including non-linear effects. It is defined by
  \begin{equation}
    \boldsymbol{u} = u, \quad\quad \boldsymbol{f}(\boldsymbol{u}) = \frac{u^2}{2},
    \label{eqn::burgers}
  \end{equation}
  with $u$ representing the flow velocity.

  \item The \emph{shallow-water equations}, also considered without viscosity, which compose a system describing a fluid flow in the regime where the vertical length scale is negligible with respect to the horizontal one. The associated state and flux are, respectively
  \begin{equation}
    \boldsymbol{u} = \left[
    \begin{array}{c}
      h\\
      hu
    \end{array}
    \right], \quad\quad \boldsymbol{f}(\boldsymbol{u}) = \left[
    \begin{array}{c}
      hu\\
      hu^2 + \frac{1}{2}gh^2
    \end{array}
    \right],
    \label{eqn::shallowWater}
  \end{equation}
  with $g$ being the gravitational constant, and $h$ the height of the fluid column.

  \item The \emph{Euler equations}, another system governing compressible fluid flow, with
  \begin{equation}
    \boldsymbol{u} = \left[
    \begin{array}{c}
      \rho\\
      \rho u\\
      E
    \end{array}
    \right], \quad\quad \boldsymbol{f}(\boldsymbol{u}) = \left[
    \begin{array}{c}
      \rho u\\
      \rho u^2 + p\\
      (E+p)u
    \end{array}
    \right]
    \label{eqn::Euler}
  \end{equation}
  where $\rho$ is the flow density, $E$ its total internal energy, and $p$ its pressure. In our case, this system is closed by considering an ideal monoatomic gas, which gives the following relationship between pressure and energy:
  \begin{equation}
    p = (\gamma-1) \left(E-\frac{1}{2}\rho u^2\right),\quad\text{with}\quad \gamma=\frac{5}{3}.
  \end{equation}
\end{itemize}
These three problems are presented in detail in \cite{leveque}, in Chap.~3.2, Chap.~5.4, and Chap.~5.1 respectively.

In the following, we provide a description of the numerical schemes employed in our experiments in order to recover their approximate solution.


\subsection{Space discretisation}
To simplify our notation, in this section we consider the state of the conservation law as being a scalar, unless otherwise specified, although the treatment described here can be seamlessly extended to the systems listed above.

Equation \eqref{eqn::consLaw} is discretised via a method of lines approach, using finite volumes in space. The spatial domain $\left[0,L\right]$ is subdivided into $N_x$ cells of uniform length $\Delta x=L/N_x$. The semi-discretised unknown is approximated by a vector $\boldsymbol{u}(t)\in\mathbb{R}^{N_x}\approx u(x,t)$, the $i$-th component of which represents the cell-average
\begin{equation}
u_i(t)\approx\frac{1}{\Delta x}\int_{x_{i-\frac{1}{2}}}^{x_{i+\frac{1}{2}}}u(x,t)\,dx,\quad\quad i=0,\dots, N_x-1.
\end{equation}
Here, $x_{i\pm\frac{1}{2}}$ identifies the coordinate of the right (respectively, left) interface of the $i$-th cell,
\begin{equation}
  x_{i\pm\frac{1}{2}}=\left(i+\frac{1\pm1}{2}\right)\Delta x,\quad
  \ifcompileFigs
    \begin{tikzpicture}[baseline, every node/.style={scale=0.7}]
      \tikzmath{
        real \cellSize; \cellSize=1.36;
        {\draw[dashed,thick,black](-1.5*\cellSize,0)--(0,0);};
        {\draw[|-,thick,black](-\cellSize,0)--(0,0) node[midway, anchor=south, black] {$i=N_x-1$};};
        {\draw (0,0) node[anchor=north, black] {$x_{-\frac{1}{2}}=0$};};
        {\draw (0,0) node[anchor=south, black] {$\equivto{x_{N_x-\frac{1}{2}}=L}{}$};};
        {\draw[|-|,thick,black] (0,0)--(\cellSize,0) node[midway, anchor=south, black] {$i=0$};};
        {\draw (\cellSize,0) node[anchor=north, black] {$x_{\frac{1}{2}}=\Delta x$};};
        {\draw[dashed,thick,black] (\cellSize,0)--(1.5*\cellSize,0);};
      }
    \end{tikzpicture}\\
  \else
    \includegraphics[valign=c]{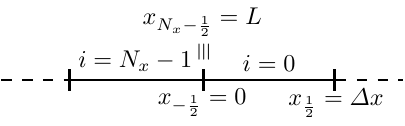}
  \fi
\end{equation}
where the right-interface of cell $N_x-1$, (on coordinate $x_{N_x-\frac{1}{2}}=L$), logically coincides with the left-interface of cell $0$, (on coordinate $x_{-\frac{1}{2}}=0$), due to the periodic boundary conditions being applied.
Integrating equation \eqref{eqn::consLaw} on cell $i$, we obtain
\begin{equation}
\frac{\partial u_i(t)}{\partial t} = - \frac{f_{i+\frac{1}{2}}(t) - f_{i-\frac{1}{2}}(t)}{\Delta x}\eqqcolon \left[\mathcal{A}(\boldsymbol{u}(t))\right]_i,
\label{eqn::spaceDisc}
\end{equation}
where $f_{i+\frac{1}{2}}(t)$ defines an approximation of the flux at the cell interface at $x_{i+\frac{1}{2}}$:
\begin{equation}
f_{i+\frac{1}{2}}(t) \approx f\left(u\left(x_{i+\frac{1}{2}},t\right)\right).
\label{eqn::fluxApprox}
\end{equation}
Typically, this approximation must be built starting from values of the numerical solution on the surrounding cells. There is, by far, no unique way to achieve this and, indeed, several schemes have been proposed in the literature, based on different methods for reconstructing the values in \eqref{eqn::fluxApprox}. Our experiments are based primarily on the WENO scheme, which is briefly introduced in \secref{sec::WENO}. In \secref{sec::LF} and \secref{sec::ROE}, we illustrate the two procedures employed to recover the fluxes in \eqref{eqn::fluxApprox} using WENO.

\subsubsection{WENO reconstruction}
\label{sec::WENO}
The \emph{Weighted Essentially Non-Oscillatory} scheme (or WENO) was originally developed by Liu, Chan and Osher in \cite{WENO}. In this section, we aim at providing a short description of the main idea behind the method, but we refer to \cite[Chap.~11.3]{jan} for a more complete treatment.

The aim of the WENO procedure is to provide a stable, high-order approximation of a certain value of interest (usually the state itself or, directly, the flux) at the cell interface. We represent these generic recovered quantities as
\begin{equation}
q^{\pm}_{i+\frac{1}{2}} \approx q\left(x_{i+\frac{1}{2}}^{\pm}\right).
\label{eqn::LRrec}
\end{equation}
Since discontinuities might arise in the solution of hyperbolic equations, the approximations on the left ($-$) and right ($+$) of an interface can differ in general, hence, the different superscripts in \eqref{eqn::LRrec}.
These quantities are approximated via polynomial reconstruction, starting from the cell-values $q_i$ in a stencil containing the interface. However, as pointed out in \secref{sec::intro}, solutions to \eqref{eqn::consLaw} can develop shocks: if the stencil used for the polynomial reconstruction happens to contain a sharp discontinuity, then the reconstructed polynomial can show oscillations due to Gibbs phenomenon, and the approximation recovered might be of poor quality. In order to counteract this, the WENO method collects approximations built using a number of different stencils, and opportunely combines them via a convex linear combination, with the weight of each approximation considered depending on an estimate of the smoothness of the particular reconstruction.

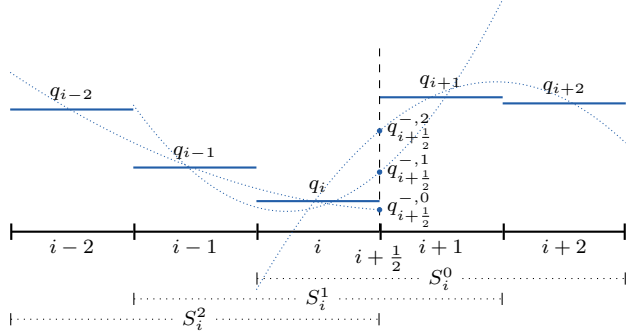
\begin{figure}[ht]
\centering
\ifcompileFigs
  \begin{tikzpicture}[scale=0.81, every node/.style={scale=0.9}]
    \def\u{{2,1.05,0.5,2.2,2.1}};     
    \tikzmath{
      int \k,\i,\temp;\k=2;
      real \cellSize; \cellSize=2;
      for \i in {0,...,2*\k-1}{
        real \uval; \uval = {\u[\i]};
        if \i == \k then {
          { \draw[|-,thick,black] ({(\i-\k-1)*\cellSize},0)--({(\i-\k)*\cellSize},0) node[midway, anchor=north, black] {$i$}; }; 
          { \draw[thick, YlGnBu-J] ({(\i-\k-1)*\cellSize},\uval)--({(\i-\k)*\cellSize},\uval)  node[midway, anchor=south, black] {$q_i$}; };
        };
        if \i < \k then {
          \temp = \k - \i;
          { \draw[|-,thick,black] ({(\i-\k-1)*\cellSize},0)--({(\i-\k)*\cellSize},0) node[midway, anchor=north, black] {$i-{\pgfmathprintnumber{\temp}}$}; };
          { \draw[thick, YlGnBu-J] ({(\i-\k-1)*\cellSize},\uval)--({(\i-\k)*\cellSize},\uval)  node[midway, anchor=south, black] {$q_{i-{\pgfmathprintnumber{\temp}}}$}; };
        };
        if \i > \k then {
          \temp = \i - \k;
          { \draw[|-,thick,black] ({(\i-\k-1)*\cellSize},0)--({(\i-\k)*\cellSize},0) node[midway, anchor=north, black] {$i+{\pgfmathprintnumber{\temp}}$}; };
          { \draw[thick, YlGnBu-J] ({(\i-\k-1)*\cellSize},\uval)--({(\i-\k)*\cellSize},\uval)  node[midway, anchor=south, black] {$q_{i+{\pgfmathprintnumber{\temp}}}$}; };
        };
      };
      \i= 2*\k;
      \temp = \i - \k;
      \uval = \u[\i];
      { \draw[|-|,thick,black] ({(\i-\k-1)*\cellSize},0)--({(\i-\k)*\cellSize},0) node[midway, anchor=north, black] {$i+{\pgfmathprintnumber{\temp}}$}; };
      { \draw[thick,YlGnBu-J] ({(\i-\k-1)*\cellSize},\uval)--({(\i-\k)*\cellSize},\uval) node[midway, anchor=south, black] {$q_{i+{\pgfmathprintnumber{\temp}}}$}; };
      { \draw[dashed, black] (0,3) -- (0,-0.1) node[anchor=north, black] {$i+\frac{1}{2}$}; };
      for \i in {{-\k-1},...,-1}{
        \temp = \i+\k+1;
        { \filldraw[|-|,dotted, black] ({(-\i-2*\k)*\cellSize},{(-\k-3.3-\i)/3})--({(-\i-\k+1)*\cellSize},{(-\k-3.3-\i)/3}) node[midway,black,fill=white]{$S_i^{\pgfmathprintnumber{\temp}}$}; };
      };
      function poly(\x,\uz,\uo,\ut) {
        real \a,\b,\c;
        \a = (   \ut - 2*\uo +    \uz ) / 2;
        \b =   - \ut + 3*\uo -  2*\uz;
        \c = ( 2*\ut - 7*\uo + 11*\uz ) / 6;
        return \a*\x*\x + \b*\x + \c;
      };
      for \i in {0,...,2}{
        real \uz,\uo,\ut;
        \uz = {\u[\i]}; \uo = {\u[\i+1]}; \ut = {\u[\i+2]};
        { \draw[densely dotted,smooth,YlGnBu-J,samples=100,domain=0:{\k+1}] plot({\cellSize*(\x+\i-3)},{poly(\x,\uz,\uo,\ut)}); };
        { \filldraw[YlGnBu-J] (0,{poly(\k+1-\i,\uz,\uo,\ut)}) circle (1pt) node[anchor=west,black]{$q_{i+\frac{1}{2}}^{-,\pgfmathprintnumber{\i}}$};  };
      };
    }
  \end{tikzpicture}
\else
  \includegraphics{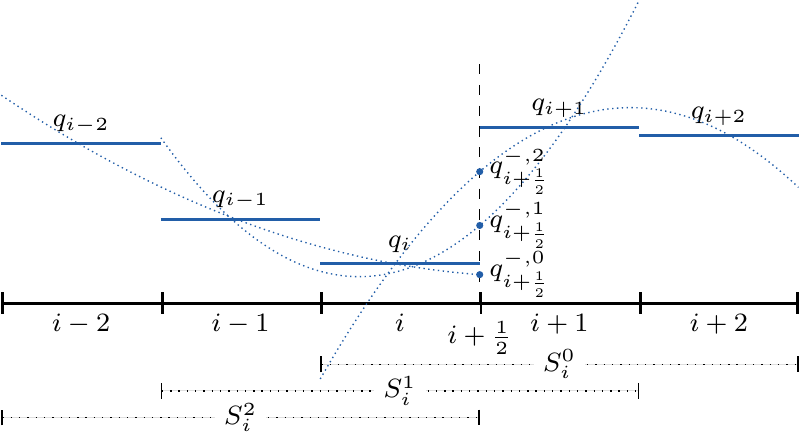}
\fi
\caption[Schematic of stencils for polynomial reconstruction in WENO]{Choice of stencils $S_i^r$ for the polynomial reconstructions of degree $k=2$ of the value left of the interface, $q^-_{i+\frac{1}{2}}$, starting from the cell values $q_{i-k},\dots,q_{i+k}$ (flat blue lines). Different stencils produce different approximating polynomials $\tilde{q}^r_i(x)$ (dotted blue lines) which, in turn, give different values at the interface (blue dots)}
\label{fig::WENO}
\end{figure}

In more detail, in order to recover an approximation of the value to the left of the interface, $q^-_{i+\frac{1}{2}}$, using a polynomial $\tilde{q}^r_i(x)$ of degree $k$, one can choose between $k+1$ different stencils each containing $k+1$ cells (see \figref{fig::WENO}).
Each of these approximations can be expressed as a linear combination of the cell-averages,
\begin{equation}
q^{-,r}_{i+\frac{1}{2}} = \tilde{q}^r_i\left(x_{i+\frac{1}{2}}\right) = \sum_{j\in S_i^r} c_{r,j-(i+r-k)}\, q_j,
\label{eqn::C}
\end{equation}
where here (and in the remainder of this section) $r$ ranges from $0$ to $k$; the stencil $S_i^r$ contains the cells with indeces $S_i^r\coloneqq\left\{i+r-k,\dots,i+r\right\}$. The quality of the $r$-th reconstruction $q^{-,r}_{i+\frac{1}{2}}$ on cell $i$ is measured via a \emph{smoothness indicator}, $\beta^r_i$, based on the (scaled) $L^2$ norm of the derivatives of the reconstructed polynomial $\tilde{q}^r_i(x)$ on $S_i^r$. This can be expressed as a bilinear form,
\begin{equation}
\beta^r_i = \sum_{l=1}^{k}\int_{x_{i-\frac{1}{2}}}^{x_{i+\frac{1}{2}}}\Delta x^{2l-1}\Bigg(\frac{\partial^l \tilde{q}^r_i(x)}{\partial x^l}\Bigg)^2dx = (\boldsymbol{q}^r_i)^T B_r\,\boldsymbol{q}^r_i,
\label{eqn::B}
\end{equation}
represented by $B_r\in\mathbb{R}^{k\times k}$, with $\boldsymbol{q}^r_i\in\mathbb{R}^{k}$ being the vector collecting the values of $q_j$ on the stencil considered, $j\in S_i^r$.
These smoothness indicators contribute to the definition of the weights $\omega^r_i$ as follows:
\begin{equation}
\omega^r_i = \frac{\alpha^r_i}{\sum_{r=0}^k\alpha^r_i},\quad\text{and}\quad\alpha^r_i=\frac{d_r}{(\varepsilon+\beta^r_i)^2},
\label{eqn::omega}
\end{equation}
where $\varepsilon$ is a small parameter to prevent division by $0$, (usually $\varepsilon=10^{-6}$). The coefficients $d_r$ are defined in such a way that, in the smooth case, it holds
\begin{equation}
  q_{i+\frac{1}{2}} = \sum_{r=0}^{k}d_r q_{i+\frac{1}{2}}^{r} = q\left(x_{i+\frac{1}{2}}\right) + \mathcal{O}(\Delta x^s),
\label{eqn::D}
\end{equation}
with $s=2k+1$. This way, the choice \eqref{eqn::omega} of the weights $\omega^r_i$ ensures that the recovered reconstruction is a high-order approximation of the underlying quantity.
Finally, the actual reconstruction is computed as a convex combination of all the different approximations in \eqref{eqn::C},
\begin{equation}
q^-_{i+\frac{1}{2}} = \sum_{r=0}^k\omega^r_i\,q^{-,r}_{i+\frac{1}{2}}.
\end{equation}
The specific values of the parameters $c_{r,j}$, $d_r$, and $B_r$ described above depend directly on the degree $k$ of the polynomials employed. For example, for $k=2$, we have:
\begin{gather}
\left[c_{r,j}\right] = \frac{1}{6}\left[ 
\begin{array}{rrr}
   2 &  5 & -1\\
  -1 &  5 &  2\\
   2 & -7 & 11
\end{array}\right],\quad \quad \left[d_r\right] = \frac{1}{10}\left[
\begin{array}{r}
  3\\
  6\\
  1
\end{array}\right],\nonumber\\
B_0 = \frac{1}{6}\left[
\begin{array}{rrr}
   20 & -31 &  11 \\
  -31 &  50 & -19 \\
   11 & -19 &   8 
\end{array}\right],\; B_1 = \frac{1}{6}\left[
\begin{array}{rrr}
    8 & -13 &   5  \\
  -13 &  26 & -13  \\
    5 & -13 &   8  
\end{array}\right],\nonumber\\
B_2 = \frac{1}{6}\left[
\begin{array}{rrr}
    8 & -19 &  11  \\
  -19 &  50 & -31  \\
   11 & -31 &  20  
\end{array}\right].
\end{gather}

An analogous procedure is used to recover the value right of the interface, $q^+_{i+\frac{1}{2}}$. We also point out that, if piecewise constant polynomials are chosen ($k=0$), then the reconstruction becomes trivial: $q^-_{i+\frac{1}{2}} = q_i$, and $q^+_{i+\frac{1}{2}} = q_{i+1}$.

When dealing with a system of conservation laws, additional care is necessary to recover a proper reconstruction of the state. A na\"ive extension from the scalar case would consider reconstructing each component of the state independently, but this has been shown to cause some unphysical oscillations in the solution \cite{WENOcharRec}. To limit this effect, a more stable approach consists, rather, of independently reconstructing the \emph{characteristic} variables of the system, since these are more readily associated with the information carried by the characteristics \cite[Chap.~11.3.4]{jan}. This comes at an additional cost, as it requires a local decomposition of the state in each of the cells considered for reconstruction, on a given reference state: this procedure is described in \cite[Procedure~2.8]{WENO2}, where an averaged state at the interface is taken as a reference. Even though we broadly follow these guidelines in our implementation, we decide for simplicity to decompose the values $\boldsymbol{u}_{i-k},\dots,\boldsymbol{u}_{i+k}$ using the central cell value $\boldsymbol{u}_{i}$ as a reference state, in order to recover $\boldsymbol{u}^{-}_{i+\frac{1}{2}}$ and $\boldsymbol{u}^{+}_{i-\frac{1}{2}}$.

With the WENO procedure available to reconstruct the states to the left and right of each interface, we proceed to approximate the numerical flux \eqref{eqn::fluxApprox} in two different ways. These are described next.

\subsubsection{Lax-Friedrichs flux}
\label{sec::LF}
We first consider the Lax-Friedrichs definition for the numerical flux \cite[Chap.~2.3.1]{WENO2}. For a system of conservation laws, this is given by
\begin{equation}
  \begin{split}
    \boldsymbol{f}^{LF}_{i+\frac{1}{2}}& = \frac{1}{2}\left(\boldsymbol{f}\left(\boldsymbol{u}_{i+\frac{1}{2}}^+\right)+\boldsymbol{f}\left(\boldsymbol{u}_{i+\frac{1}{2}}^-\right)\right)\\
                                       & + \frac{1}{2}\alpha_{LF} \left(\boldsymbol{u}_{i+\frac{1}{2}}^+ - \boldsymbol{u}_{i+\frac{1}{2}}^-\right).
  \end{split}
  \label{eqn::LFflux}
\end{equation}
This flux is one of the simplest to prescribe; however, as a downside, it typically produces smeared out numerical solutions, since it adds artificial diffusion to the system. The amount of numerical diffusion added depends directly on the parameter $\alpha_{LF}$: for our implementation, we choose the (somewhat loose) value of
\begin{equation}
  \alpha_{LF} = \max_{i,k}\left\{\left|\lambda_k\left(\boldsymbol{u}^{\pm}_{i+\frac{1}{2}}\right)\right|\right\},
  \quad\quad\begin{array}{ll} 
    i=& 0,\dots,N_x-1\\
    k=& 0,\dots,D-1
    \end{array},
\end{equation}
where the $\lambda_k$'s are the $D$ eigenvalues (if the system is composed of $D$ equations) of the Jacobian of the flux:
\begin{equation}
  J_{\boldsymbol{f}}\left(\boldsymbol{u}\right)\boldsymbol{r}_k=\lambda_{k}\boldsymbol{r}_k,\quad\text{with}\quad\left[J_{\boldsymbol{f}}\left(\boldsymbol{u}\right)\right]_{m,n}=\frac{\partial\left[\boldsymbol{f}(\boldsymbol{u})\right]_m}{\partial\left[\boldsymbol{u}\right]_n},
  \label{eqn::eigs}
\end{equation}
while the $\boldsymbol{r}_k$'s represent the corresponding eigenvectors.


\subsubsection{Roe flux}
\label{sec::ROE}
The second scheme investigated in our experiments uses Roe's approximate Riemann solver \cite[Chap.~14.2]{leveque} in order to recover the numerical flux. For each interface, this solver targets a linearisation of the Riemann problem defined by \eqref{eqn::consLaw} and the left and right values computed using the WENO procedure. This problem is written as:
\begin{equation}
  \boldsymbol{u}_t + J_{\boldsymbol{f}}(\hat{\boldsymbol{u}})\,\boldsymbol{u}_x = 0,
  \label{eqn::linearConsLaw}
\end{equation}
where $\hat{\boldsymbol{u}}$ is an opportunely defined \emph{Roe-averaged} state.
All relevant variables in this section are to be intended as evaluated at an interface: we drop the subscripts $i+1/2$ for ease of notation.
The solution to the linear Riemann problem \eqref{eqn::linearConsLaw}, as well as the associated flux at the interface, are both easily expressed in terms of the eigenvalues and eigenvectors of the Jacobian. These are given by
\begin{equation}
  \begin{array}{l}
    \lambda_0^{S} = u-c_{S}\\
    \lambda_1^{S} = u+c_{S}
  \end{array} \quad\text{and}\quad
  \begin{array}{l}
    \boldsymbol{r}_0^{S} = \left[\begin{array}{cc}1, & \lambda_0^{S}\end{array}\right]^T\\
    \boldsymbol{r}_1^{S} = \left[\begin{array}{cc}1, & \lambda_1^{S}\end{array}\right]^T
  \end{array}
\end{equation}
for the shallow-water system, with \emph{speed of sound} $c_{S}=\sqrt{gh}$, and by
\begin{equation}
  \begin{array}{l}
    \lambda_0^{E} = u-c_E\\
    \lambda_1^{E} = u    \\
    \lambda_2^{E} = u+c_E
  \end{array} \quad\text{and}\quad
  \begin{array}{l}
    \boldsymbol{r}_0^{E} = \left[\begin{array}{ccc}1, & \lambda_0^{E}, & H-uc_E \end{array}\right]^T\\
    \boldsymbol{r}_1^{E} = \left[\begin{array}{ccc}1, & \lambda_1^{E}, & u^2/2  \end{array}\right]^T\\
    \boldsymbol{r}_2^{E} = \left[\begin{array}{ccc}1, & \lambda_2^{E}, & H+uc_E \end{array}\right]^T
  \end{array}
\end{equation}
for the Euler equations, with speed of sound $c_{E}=\sqrt{(\gamma-1)(H-u^2/2)}$. Here, $H=(E+p)/\rho$ is the \emph{specific enthalpy}. For Burgers' equation, the only eigenvalue is given by $\lambda_0^B = u$ itself.
Since we are only interested in the eigenvalues and eigenvectors of the Jacobian $J_{\boldsymbol{f}}(\hat{\boldsymbol{u}})$ in \eqref{eqn::linearConsLaw}, (namely $\hat{\lambda}_k$ and $\hat{\boldsymbol{r}}_k$), it suffices to define the following Roe-averaged quantities: an average velocity for the Burgers' equation, given by
\begin{equation}
  \hat{u}=\frac{u^+ + u^-}{2};
\end{equation}
average height and velocity for the Shallow-water equation, prescribed as
\begin{equation}
  \hat{h} = \frac{h^++h^-}{2} \quad\text{and}\quad 
  \hat{u} = \frac{u^+\sqrt{h^+} + u^-\sqrt{h^-}}{\sqrt{h^+}+\sqrt{h^-}};
\end{equation}
finally, average velocity and specific enthalpy for the Euler equations:
\begin{equation}
  \begin{split}
    \hat{u}&=\frac{u^+\sqrt{\rho^+} + u^-\sqrt{\rho^-}}{\sqrt{\rho^+}+\sqrt{\rho^-}} \quad\text{and}\quad\\
    \hat{H}&=\frac{(E^++p^+)\sqrt{\rho^+} + (E^-+p^-)\sqrt{\rho^-}}{\sqrt{\rho^+}+\sqrt{\rho^-}}.
  \end{split}
\end{equation}
Since the target problem \eqref{eqn::linearConsLaw} is linear, this procedure is known to provide a non-entropic weak solution when transonic rarefaction waves arise \cite[Chap.~14.2.2]{leveque}. To counteract this, we also apply the \emph{entropy fix} proposed by Harten and Hyman in \cite{HHentropyFix}.
In the general case of a system of $D$ conservation laws, this gives rise to the following formula for the numerical flux at the interface:
\begin{equation}
\boldsymbol{f}^{R} = \frac{1}{2}\left( \boldsymbol{f}\left(\boldsymbol{u}^-\right) + \boldsymbol{f}\left(\boldsymbol{u}^+\right)\right)\, -\,\sum_{k=0}^{D-1} q_{H}(\hat{\lambda}_k)\,\alpha_k\,\hat{\boldsymbol{r}}_k.
\label{eqn::roeflux}
\end{equation}
Here, $\alpha_k$ is the \emph{shock strength}, that is the jump in the $k$-th characteristic variable between the states left and right of the interface:
$\alpha_k = \left(\boldsymbol{u}^+ - \boldsymbol{u}^-\right)\cdot\hat{\boldsymbol{r}}_k$.
Also, $q_{H}(\lambda_k)$ defines the actual entropy fix:
\begin{equation}
  q_{H}(\hat{\lambda}_k) = \left\{
  \begin{array}{ll}
    \frac{1}{2}\left(\frac{\hat{\lambda}_k}{\delta_k}+\delta_k\right) &\quad\quad \text{if}\quad  |\hat{\lambda}_k|<\delta_k\\
    |\hat{\lambda}_k| &\quad\quad \text{else} 
  \end{array}\right.,
\end{equation}
where $\delta_k=\max\left\{0,\hat{\lambda}_k-\lambda_k^-,\lambda_k^+-\hat{\lambda}_k\right\}$, and $\lambda_k^{\pm}$ are the eigenvalues of $J_{\boldsymbol{f}}(\boldsymbol{u}^\pm)$. 


\paragraph*{}
Substituting \eqref{eqn::LFflux} or \eqref{eqn::roeflux} into \eqref{eqn::spaceDisc} identifies the right-hand side $\mathcal{A}(\boldsymbol{u})$ of our semi-discretised equation. The details regarding the discretisation in time, as well as the definition of the time-steppers employed, are presented next.



\subsection{Time discretisation}
The time domain is also discretised using a uniform grid of $N_t$ nodes, with time step $\Delta t = T/N_t$. The unknowns at each instant $t_n=n\Delta t$, $n=0,\dots,N_t-1$ are approximated by a vector denoted as $\boldsymbol{u}^n\approx\boldsymbol{u}(t_n)$.

Since we aim at employing high-order spatial reconstructions, it is sensible to request that our temporal discretisation matches this accuracy. The schemes chosen belong to the family of \emph{Strong Stability-Preserving Runge Kutta} methods (SSPRK). As the name suggests, these have the remarkable property of being able to preserve strong stability and, in particular, to be \emph{Total Variation Diminishing}, even while achieving high-order of accuracy; for this reason, they are often employed in conjunction with high-order spatial discretisation of hyperbolic equations. We refer to \cite{SSPRK} for a complete review of these schemes, and we only report here the definition of the methods used in our experiments, which vary in the order of accuracy $d$ recovered. They are: a third-order scheme (SSPRK3), whose stepping procedure applied to \eqref{eqn::spaceDisc} is given by
\begin{equation}
\begin{split}
  \boldsymbol{u}^{n,(1)} &=                                            \boldsymbol{u}^{n}     +            \Delta t\, \mathcal{A}\left(\boldsymbol{u}^{n}    \right)\\
  \boldsymbol{u}^{n,(2)} &= \frac{3}{4}\boldsymbol{u}^{n} + \frac{1}{4}\boldsymbol{u}^{n,(1)} + \frac{1}{4}\Delta t\, \mathcal{A}\left(\boldsymbol{u}^{n,(1)}\right)\\
  \boldsymbol{u}^{n+1}   &= \frac{1}{3}\boldsymbol{u}^{n} + \frac{2}{3}\boldsymbol{u}^{n,(2)} + \frac{2}{3}\Delta t\, \mathcal{A}\left(\boldsymbol{u}^{n,(2)}\right);
\end{split}
\label{eqn::SSPRK3}
\end{equation}
a second-order scheme (SSPRK2), defined by
\begin{equation}
\begin{split}
  \boldsymbol{u}^{n,(1)} &=            \boldsymbol{u}^{n}                                     +            \Delta t\, \mathcal{A}\left(\boldsymbol{u}^{n}    \right)\\
  \boldsymbol{u}^{n+1}   &= \frac{1}{2}\boldsymbol{u}^{n} + \frac{1}{2}\boldsymbol{u}^{n,(1)} + \frac{1}{2}\Delta t\, \mathcal{A}\left(\boldsymbol{u}^{n,(1)}\right);
\end{split}
\label{eqn::SSPRK2}
\end{equation}
and, finally, the first-order scheme, which reduces to the well-known Forward Euler (FE) method,
\begin{equation}
  \boldsymbol{u}^{n+1} = \boldsymbol{u}^{n} + \Delta t\, \mathcal{A}\left(\boldsymbol{u}^{n}\right).
  \label{eqn::SSPRK1}
\end{equation}

\section{MGRIT}
\label{sec::MGRIT}
The Multigrid Reduction In Time method can be interpreted as a multigrid scheme with the coarsening procedure acting along the time domain. First introduced in \cite{MGRIToriginal}, it has quickly become a mature algorithm for time-parallelisation. In this section, we provide a brief description of the scheme.

\subsection{Method overview}
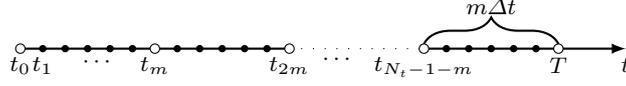
\begin{figure}
\centering
\ifcompileFigs
  \begin{tikzpicture}[scale=0.59]
  \tikzmath{
    \nf = 6;
    \Np = 2;
    \DT = 3;
    \rg = 3;
    \rf = 2;
    for \i in {0,...,\Np-1}{
      {\draw[thick] ( \i*\DT, 0 ) -- ( \i*\DT+\DT, 0 );};
      {\draw[fill=white] ( \i*\DT, 0 ) circle (\rg pt);};
      for \j in {1,...,\nf-1}{
        {\draw[fill=black] ({\i*\DT+\j*\DT/(\nf)}, 0) circle (\rf pt) ;};
      };
    };
    {\draw[fill=black] (0, 0) node[anchor=north]{$t_{0}$};};
    {\draw[fill=black] ({\DT/(\nf)}, 0) node[anchor=north]{$t_{1}$};};
    {\draw[fill=black] ({(1+\nf)/2*\DT/(\nf)}, 0) node[anchor=north]{$\cdots$};};
    {\draw[fill=black] ({\DT}, 0) node[anchor=north]{$t_m$};};
    {\draw[loosely dotted] ( \Np*\DT, 0 ) -- ( \Np*\DT+\DT, 0 ) node[anchor=north,pos=.35]{$\ldots$};};
    {\draw[fill=white] ( \Np*\DT, 0 ) circle (\rg pt) node[anchor=north]{$t_{\pgfmathprintnumber{\Np}m}$};};
    {\draw[thick,->,>=latex] ( \Np*\DT+\DT, 0 ) -- ( \Np*\DT+2*\DT+\DT/2, 0 ) node[anchor=north]{$t$};};
    {\draw[fill=white] ( \Np*\DT+\DT, 0 ) circle (\rg pt) node[anchor=north]{$t_{N_t-1-m}$};};
    {\draw[fill=white] ( \Np*\DT+2*\DT, 0 ) circle (\rg pt) node[anchor=north]{$T$};};
    {\draw [decorate,decoration={brace,amplitude=10pt},thick,yshift=3 pt]
    ( \Np*\DT+\DT, 0 ) -- ( \Np*\DT+2*\DT, 0 ) node [black,midway,anchor=south,yshift =8 pt] 
    {$m\Delta t$};};
    for \j in {1,...,\nf-1}{
      {\draw[fill=black] ({\Np*\DT+\DT+\j*\DT/(\nf)}, 0) circle (\rf pt) ;};
    };
  }
  \end{tikzpicture}
\else
  \includegraphics{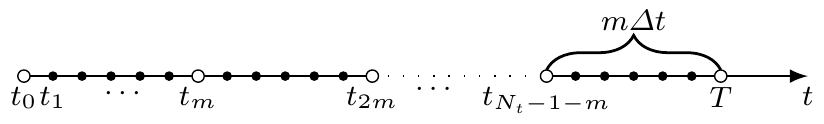}
\fi
\caption[Sketch of the partitioning of the time discretisation]{Sketch of the partitioning of the time discretisation: white dots represent coarse nodes, black dots are fine nodes. The last time chunk is also highlighted}
\label{fig::timeGridsMGRIT}
\end{figure}

As the name suggests, MGRIT is, at its root, a multigrid reduction technique applied to the monolithic system arising from the space-time discretisation of a PDE such as \eqref{eqn::consLaw}. In the constant-coefficient, linear, homogeneous case with a single-step time integrator, such a system takes the form of a block bi-diagonal, block Toeplitz matrix:
\begin{equation}
A\boldsymbol{u} = \boldsymbol{g} \Leftrightarrow \left[
\begin{array}{cccc}
  I      &        &        &\\
  -\Phi  & I      &        &\\
         & \ddots & \ddots &\\
         &        & -\Phi  & I
\end{array}\right]\left[
\begin{array}{c}
  \boldsymbol{u}^0\\ 
  \boldsymbol{u}^1\\
  \vdots          \\
  \boldsymbol{u}^{N_t-1}
\end{array}\right] = \left[
\begin{array}{c}
  \bar{\boldsymbol{u}}^0\\ 
  \boldsymbol{0}\\
  \vdots        \\
  \boldsymbol{0}
\end{array}\right],
\label{eqn::monolithicSystem}
\end{equation}
where $\boldsymbol{u}\in\mathbb{R}^{N_xN_t}$ is the vector containing the values of the discrete solution at each node in the space-time grid, and $\bar{\boldsymbol{u}}^0$ contains the discretisation of the initial condition; since we do not consider forcing terms, the rest of $\boldsymbol{g}$ is filled with zeros. The \emph{fine integrator} $\Phi$ represents the action of the time-stepper, so that, for each instant, $n$, we have $\boldsymbol{u}^{n+1} = \Phi\boldsymbol{u}^{n}$. Following the multigrid philosophy, the nodes composing the temporal discretisation are split into two sets, respectively, of \emph{coarse} (denoted with $C$) and \emph{fine} (denoted with $F$) nodes. In our case, we simply pick the coarse nodes to be spaced every $m$ nodes, with $m$ being the \emph{coarsening factor} (see \figref{fig::timeGridsMGRIT}). The coarse nodes thus subdivide the time domain into $(N_t-1)/m$ different \emph{time chunks}, each of size $m\Delta t$.
The variables in the monolithic space-time discretisation $\boldsymbol{u}$ and the coefficients in $A$ can also be rearranged accordingly, so that the matrix can be factorised in the following way:
\begin{equation}
\begin{split}
A &= \left[
  \begin{array}{cc}
    A_{FF} & A_{FC} \\
    A_{CF} & A_{CC}
  \end{array}
\right] \\
  &= \left[
  \begin{array}{cc}
    I_F               & 0 \\
    A_{CF}A_{FF}^{-1} & I_C
  \end{array}
\right] \left[
  \begin{array}{cc}
    A_{FF} & 0 \\
    0      & S
  \end{array}
\right]\left[
  \begin{array}{cc}
    I_F & A_{FF}^{-1}A_{FC} \\
    0   & I_C
  \end{array}
\right],
\end{split}
\label{eqn:AfactMGRIT}
\end{equation}
where $I_C$ and $I_F$ are identity matrices of appropriate sizes. 
This factorisation allows us to separate the solution over the $F$ and the $C$ nodes. It can be easily shown that the $A_{FF}$ sub-matrix in \eqref{eqn:AfactMGRIT} presents a block-diagonal structure, which implies that systems involving it can be readily solved in parallel:
\begin{equation}
  A_{FF} =  \left[
  \begin{array}{ccc}
    B_\Phi &        & \\
           & \ddots & \\
           &        & B_\Phi
  \end{array}\right],\quad B_\Phi = \left[
  \begin{array}{cccc}
    I     &        &        &\\
    -\Phi & I      &        &\\
          & \ddots & \ddots &\\
          &        & -\Phi  & I
  \end{array}\right].
  \label{eqn::AFF}
\end{equation}
Here, we note that $A_{FF}$ is, overall an $N_x(N_t-1)(m-1)/m \times N_x(N_t-1)(m-1)/m$ matrix, written above as an $(N_t-1)/m \times (N_t-1)/m$ block diagonal system, with blocks $B_\Phi$ of size $(m-1)N_x \times (m-1)N_x$. Each solve with $B_\Phi$ involves ``stepping'' forward an approximate solution on the spatial mesh through $m-2$ applications of $\Phi$; each of these solves is completely independent and can be performed in parallel.
The challenge lies rather in finding the solution to the system associated with the Schur complement: 
\begin{equation}
  S=A_{CC}-A_{CF}A_{FF}^{-1}A_{FC} =  \left[
  \begin{array}{cccc}
    I        &        &          &\\
    -\Phi^m  & I      &          &\\
             & \ddots & \ddots   &\\
             &        & -\Phi^m  & I
  \end{array}\right].
  \label{eqn::schurComplement}
\end{equation}
Attempting to solve this directly is equivalent to applying forward block-substitution to \eqref{eqn::monolithicSystem}, (\emph{i.e.}, time-stepping over all temporal nodes) and would, thus, nullify any advantage gained from parallelisation. Rather, in the scope of the MGRIT algorithm, we resort to solving a modified system, obtained by substituting another operator $\Psi\approx\Phi^m$ in \eqref{eqn::schurComplement}, giving
\begin{equation}
  S_{\Psi} =  \left[
  \begin{array}{cccc}
    I      &        &        &\\
    -\Psi  & I      &        &\\
           & \ddots & \ddots &\\
           &        & -\Psi  & I
  \end{array}\right]\approx S.
  \label{eqn::approxSchurComplement}
\end{equation}
The structure of this system is equivalent to that of \eqref{eqn::monolithicSystem}, so that $\Psi$ can be interpreted as yet another time-stepping routine, which acts only on the coarse nodes: we call this operator the \emph{coarse integrator}. Notice that the splitting into fine and coarse nodes can be applied in a recursive fashion, further extracting a hierarchy of $N_l$ grids, together with their corresponding integrators $\Phi_{(l)}$, in a true \emph{multi}-grid spirit.

The hope is that, by opportunely alternating between solving for the fine nodes (inverting $A_{FF}$) and for the coarse nodes (inverting $S_{\Psi}$) and iterating, we can quickly converge to a suitable approximation of the solution of the original system \eqref{eqn::monolithicSystem}.

\subsection{The algorithm}
In more detail, an iteration of the MGRIT algorithm consists of the following building blocks:

\paragraph{Relaxation} Update the solution at the fine nodes of the current level, given a guess at the coarse nodes. This involves solving a system in the form
\begin{equation}
  A_{FF}^{(l)}\boldsymbol{u}_{(l)}^F=\boldsymbol{g}_{(l)}^F - A_{FC}^{(l)}\boldsymbol{u}_{(l)}^C,
  \label{eqn::relaxation}
\end{equation}
where $A_{FF}^{(l)}$ has the same structure as in \eqref{eqn::AFF}, but whose subdiagonal blocks contain $\Phi_{(l)}$, while $\boldsymbol{u}_{(l)}^F$ and $\boldsymbol{u}_{(l)}^C$, respectively, represent the grouping of the unknowns at the fine and coarse nodes of the current level, $l$. The solution is, hence, updated by time-stepping using $\Phi_{(l)}$, starting from the given values at the coarse nodes. Here lies the parallel part of the algorithm, since the time-stepping procedure can be carried out independently on each chunk. In our work, we consider two types of relaxation: \emph{F-relaxation}, where the time-stepping covers a single time chunk, updating the fine nodes within it, and \emph{FCF-relaxation}, where the time-stepping carries on over the following coarse node (performing a \emph{C-relaxation}) and continues on the following time chunk (adding another F-relaxation).  This is an \emph{overlapping} form of relaxation that requires more work per level of the hierarchy but can be implemented with the same parallel efficiency.

\paragraph{Restriction} Transfer information from the nodes of the current level $l$ to the nodes at the coarser level $l+1$. Simple injection onto the coarse nodes is chosen as the restriction operator:
\begin{equation}
  R_{(l)}^{(l+1)}\boldsymbol{u}_{(l)} = \boldsymbol{u}_{(l)}^C.
  \label{eqn::restrictionOperator}%
\end{equation}
Since we are dealing with non-linear equations, we also implement the \emph{Full Approximation Storage} (FAS) algorithm \cite[Chap.~5.3.4]{multigrid}. This modifies the right-hand side of the system at the coarse level by adding a correction term:
\begin{equation}
  \begin{split}
    S_{\phi_{(l+1)}}\left(\boldsymbol{u}_{(l+1)}\right) & = R_{(l)}^{(l+1)}\left(\boldsymbol{g}_{(l)} - S_{\phi_{(l)}}\left(\boldsymbol{u}_{(l)}\right)\right) \\
                                                        & + \underbrace{S_{\phi_{(l+1)}}\left(R_{(l)}^{(l+1)}\boldsymbol{u}_{(l)}\right)}_{\text{FAS correction}} = \boldsymbol{g}_{(l+1)},
    \label{eqn::FAScorrection}%
  \end{split}
\end{equation}
where $S_{\Phi_{(l+1)}}$ has the same structure as in \eqref{eqn::approxSchurComplement}, but with $\Phi_{l+1}$ on the subdiagonal. We point out that the same restriction operator \eqref{eqn::restrictionOperator} is applied to both the residual and the state.
Given the particular structure of the operators involved, formula \eqref{eqn::FAScorrection} simplifies, so that the right-hand side of the coarse-level system, $\boldsymbol{g}_{(l+1)}$, reduces to
\begin{equation}
  \boldsymbol{g}_{(l+1)}^i = \boldsymbol{g}_{(l)}^{mi} + \Phi_{(l)}\left(\boldsymbol{u}_{(l)}^{mi-1}\right) - \Phi_{(l+1)}\left(\boldsymbol{u}_{(l)}^{m(i-1)}\right),
  \label{eqn::restriction::rhs}
\end{equation}
for each temporal node $i=1,\dots,(N_t-1)/m^{l+1}+1$ of the coarser level. We also need to provide an initial guess for the solution at the coarse level $\boldsymbol{u}_{(l+1)}$: this is usually chosen to be $\boldsymbol{u}_{(l)}^C$, however we have found that a better alternative is given by
\begin{equation}
  \boldsymbol{u}_{(l+1)}^i = \Phi_{(l)}\left(\boldsymbol{u}_{(l)}^{mi-1}\right) + \boldsymbol{g}_{(l)}^{mi},
  \label{eqn::restriction::u}
\end{equation}
that is, by performing the last integration step for each chunk and injecting the resulting value. This operation comes at virtually no additional cost (since the quantity in \eqref{eqn::restriction::u} is already computed as part of \eqref{eqn::restriction::rhs}), and in early experiments it was seen to improve convergence: in \figref{fig::MGRITrestrictionGuess}, we give an example of the effectiveness of \eqref{eqn::restriction::u} over simple injection.
Notice that this procedure is \emph{not} equivalent to performing an additional C-relaxation before injection, since the FAS correction \eqref{eqn::FAScorrection} is still based on the non-updated values at the coarse nodes $\boldsymbol{u}_{(l)}^{m(i-1)}$; however, this approach \emph{does} retain the fixed-point property of the FAS algorithm since \eqref{eqn::restrictionOperator} and \eqref{eqn::restriction::u} coincide when $\boldsymbol{u}_{(l)}$ is the exact solution on level $l$.

\paragraph{Coarse grid correction} Update solution at the coarsest level. This involves solving the system
\begin{equation}
  S_{\Phi_{(N_l-1)}}\boldsymbol{u}_{(N_l-1)}=\boldsymbol{g}_{(N_l-1)}.
  \label{eqn::coarseCorr}
\end{equation}
This procedure is sequential but, by choosing a cheap coarse integrator, the overhead to the algorithm can be limited, and parallel efficiency can still be gained.

\paragraph{Interpolation} Transfer information from the nodes at the coarser level $l+1$ to those at the current level $l$. \emph{Ideal interpolation} is chosen, which in our case consists of two steps: an injection from coarse to fine grid, followed by an F-relaxation on the fine level, starting from the freshly updated values at the coarse nodes. Overall, then, we have the following definition for the interpolation operator:
\begin{equation}
  \begin{split}
    &\left[I^{(l)}_{(l+1)}\left(\boldsymbol{u}_{(l+1)}\right)\right]^i = F_{(l)}^{\floor*{\frac{i}{m}},i\%m}\left(\boldsymbol{u}_{(l+1)}^{\floor*{\frac{i}{m}}}\right), \text{ with}\\
    &F_{(l)}^{k,r}(\boldsymbol{u}) =
    \begin{cases}
      \boldsymbol{u} & \text{if }r=0,\\
      \Phi_{(l)}\left(F_{(l)}^{k,r-1}(\boldsymbol{u})\right) + \boldsymbol{g}_{(l)}^{k+r} & \text{else},
    \end{cases}\\
  \end{split}
  \label{eqn::interpolation}
\end{equation}
where $\%$ represents the \emph{modulo} operator, while $\floor{*}$ the \emph{floor} operator.
Notice that, if the interpolation step is followed by a relaxation step, then the F-relaxation within interpolation becomes redundant, since the results would get overwritten by the relaxation step.  Here, we consider only pre-relaxation within the multigrid cycles, and write the interpolation step as injection followed by an F-relaxation step.

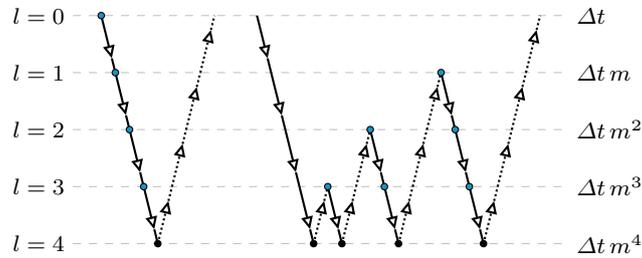
\begin{figure}[ht]
\centering
\ifcompileFigs
  \begin{tikzpicture}[scale=0.54]
  \tikzmath{
    \nlvl = 5;
    \hlvl = 1.4;
    \llvl = 0.345;
    \L    = 7*\llvl + 2*\llvl*(\nlvl-1) + \llvl*\nlvl*(\nlvl-1);
    {\draw[dashed, lightgray, line width=0.10mm] ( 0, 0        )  node[anchor=east, black]{$l=0$} -- ( {\L}, 0        ) node[anchor=west, black]{$ \Delta t$};};
    {\draw[dashed, lightgray, line width=0.10mm] ( 0, {-\hlvl} )  node[anchor=east, black]{$l=1$} -- ( {\L}, {-\hlvl} ) node[anchor=west, black]{$\Delta t\,m$};};
    for \i in {2,...,\nlvl-1}{
      {\draw[dashed, lightgray, line width=0.10mm] ( 0, {-\i*\hlvl} ) node[anchor=east, black]{$l=\pgfmathprintnumber{\i}$} -- ( {\L}, {-\i*\hlvl} ) node[anchor=west, black]{$\Delta t\,m^{\pgfmathprintnumber{\i}}$};};
    };
    function plotBaseCase(\lvl,\xshift){
      {\filldraw[black] (\xshift,{-\lvl*\hlvl}) circle (2.25pt);};
    };
    function plotRestriction(\lvl,\xshift){
      \xh = \xshift+0.75*\llvl;
      \yh = {-(\lvl+0.75)*\hlvl};
      {\draw[thick,black,-{Triangle[open, width=1.5mm]},>=Latex] (\xshift,{-\lvl*\hlvl}) -- (\xh,\yh);};
      {\draw[thick,black,] (\xh,\yh) -- ({\xshift+\llvl},{-(\lvl+1)*\hlvl});};
      {\draw[black,fill=YlGnBu-H, shorten <=2.25pt] ({\xshift},{-\lvl*\hlvl}) circle (2.25pt);};
    };
    function plotInterpolation(\lvl,\xshift){
      \xh = \xshift+0.75*\llvl;
      \yh = {-(\lvl-0.75)*\hlvl};
      {\draw[thick,black,densely dotted,-{Triangle[open, width=1.5mm]},>=latex] (\xshift,{-\lvl*\hlvl}) -- (\xh,\yh);};
      {\draw[thick,black,densely dotted,] (\xh,\yh) -- ({\xshift+\llvl},{-(\lvl-1)*\hlvl});};
    };
    function Vcycle(\lvl,\xshift) {
      if \lvl == \nlvl-1 then {
        plotBaseCase(\lvl,\xshift);
        return \xshift;
      } else {
        plotRestriction(\lvl,\xshift);
        \xshift = Vcycle(\lvl+1,\xshift+\llvl);
      };
      if \lvl > 0 then{
        plotInterpolation(\lvl,\xshift);
      };
      return \xshift+\llvl;
    };
    \xshift = 2*\llvl;
    \xshift = Vcycle(0,\xshift) + 2*\llvl;
    for \i in {0,...,\nlvl-2}{
      \xh = \xshift+0.75*\llvl;
      \yh = {-(\i+0.75)*\hlvl};
      {\draw[thick,black,-{Triangle[open, width=1.5mm]},>=latex] (\xshift,{-\i*\hlvl}) -- (\xh,\yh);};
      {\draw[thick,black,] (\xh,\yh) -- ({\xshift+\llvl},{-(\i+1)*\hlvl});};
      \xshift = \xshift + \llvl;
    };
    function Fcycle(\lvl,\xshift) {
      if \lvl == 1 then {
        return \xshift;
      } else {
        \xshift = Vcycle(\lvl-1,\xshift);
        Fcycle(\lvl-1,\xshift);
      };
    };
    Fcycle(\nlvl,\xshift);
  }
  \end{tikzpicture}
\else
  \includegraphics{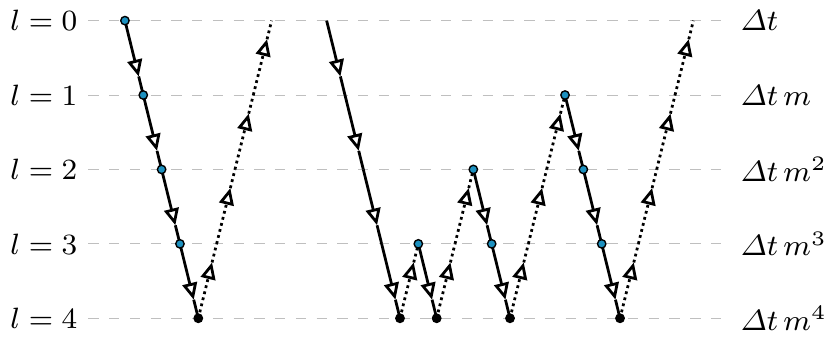}
\fi
\caption[Sketch of different multigrid cycles]{Sketch of the flow of the MGRIT algorithm if prescribing a V-cycle (left) or an F-cycle (right), on 5 levels. Restriction (full lines) is applied when moving from a level to a coarser one (dashed lines: the size of the time-step at each level $l$ is reported on the right). Relaxation (blue dots) is applied at each level, bar the bottom one, where the coarsest system is solved exactly (black dots). Interpolation (dotted lines) is applied when moving from a level to a finer one}
\label{fig::MGcycles}
\end{figure}

\paragraph*{}
The order in which the algorithm moves between discretisation levels defines the sequence in which the operations above are prescribed, and identifies the type of \emph{cycle} used in MGRIT. We experiment both with a \emph{V-cycle} and an \emph{F-cycle} \cite[Chap.~3]{multigrid}. In \algref{alg::Vcycle}, the pseudo-code outlining the (recursive) definition of a V-cycle with F-relaxation is provided; a sketch of the movement along the various levels for both V- and F-cycle is given in \figref{fig::MGcycles}.
\begin{algorithm}[ht]
  \setcounter{AlgoLine}{0}
  \uIf{$l<N_l-1$}{
    \tcc{Parallel step}
    Perform relaxation at current level and recover $\boldsymbol{u}_{(l)}^F$ by solving \eqref{eqn::relaxation}\;
    Compute residual and FAS correction, and restrict to coarser level to compute $\boldsymbol{g}_{(l+1)}$ and $\boldsymbol{u}_{(l+1)}$ following formulae \eqref{eqn::restriction::rhs} and \eqref{eqn::restriction::u}\;
    \BlankLine
    \tcc{Recursive step}
    Invoke this function at coarser level and find $\boldsymbol{u}_{(l+1)}$\;
    \BlankLine
    \tcc{Parallel step}
    Interpolate from coarser level to recover $\boldsymbol{u}_{(l)}$ applying \eqref{eqn::interpolation}\;
  }\Else{
    \tcc{Serial step - base step of recursion}
    Recover $\boldsymbol{u}_{(N_l-1)}$ by solving system \eqref{eqn::coarseCorr}\;
  }
\caption{MGRIT, V-cycle iteration}
\label{alg::Vcycle}
\end{algorithm}

\subsection{Choices for the coarse integrator}
\label{sec::integrators}
As hinted above, the key to the effective application of MGRIT lies in finding an adequate pair (or set) of fine and coarse integrators. In particular, the coarse integrator needs to do a ``good enough'' job of mimicking the action of the fine integrator, so that \eqref{eqn::approxSchurComplement} represents a valid approximation to \eqref{eqn::schurComplement}; at the same time, it should be cheap enough to compute, so that the solution to systems involving \eqref{eqn::approxSchurComplement} can be promptly recovered.
The necessity of overcoming this trade-off between cost and accuracy is common to many other time-parallelisation algorithms, and a number of approaches have been proposed to address this. The most straightforward is simple rediscretisation: the coarse integrator is none other than the fine integrator applied on a coarser time grid \cite{MGRIToriginal,MGRITNonLin}; an alternative consists in choosing integrators of varying levels of accuracy \cite{Nielsen}; yet another lies in neglecting or simplifying certain physics at the coarse level when timescale separation is possible \cite{Beth}.
These solutions focus mostly on finding a coarse integrator which is cheap to compute, but it is unclear to what extent it remains faithful to the fine integrator. Indeed, in the hyperbolic framework, simple rediscretisation has been shown to fail in many cases (see \figref{fig::MGRITfineMatch}), particularly if low-order time discretisations are employed: evidence of this is given in \cite{MGRIT_Oliver}, as well as in the results of \figref{fig::MGRITcfl}, at least for certain regimes. When sticking to simple time-steppers, then, additional care is needed to ensure stability of MGRIT. To address this, we propose an approach which reverses the aforementioned principle: we start directly from the definition of the (explicit) fine integrator, and we progressively perform some simplifications which renders it cheaper to apply, without sacrificing ``too much'' accuracy.  With this, we are able to ``restore'' MGRIT convergence for a low-order explicit time discretisation to be comparable to that of direct rediscretisation with high-order explicit time discretisations, although in neither case do we see CFL-robust convergence.

\subsubsection{Fine integrator matching}
\label{sec::fineMatching}
To illustrate our approach, we consider our model equation \eqref{eqn::burgers}. If a simple Lax-Friedrichs flux discretisation is applied in conjunction with a Forward Euler time discretisation, the resulting time-stepping formula reads
\begin{equation}
\begin{array}{rcrc}
  u_i^{n+1} =& \underbrace{\frac{u_{i+1}^n + u_{i-1}^n}{2}} & - \Delta t & \underbrace{\frac{f(u_{i+1}^n) - f(u_{i-1}^n)}{2\Delta x}}\\[15pt]
            =& Eu_i^n                                     & - \Delta t & Df_i^n,
\end{array}
\label{eqn::fineLF}
\end{equation}
with $f_i^n = f(u_i^n)$, and where we introduced the \emph{central difference} operator, $D((*)_i) = ((*)_{i+1}-(*)_{i-1})/2\Delta x$, and the \emph{average} operator $E((*)_i) = ((*)_{i+1}+(*)_{i-1})/2$.
If we were to simply rediscretise at the coarse level and apply the same scheme to a coarser grid with time step $2\Delta t$, the resulting formula would similarly read
\begin{equation}
u_i^{n+1} = Eu_i^n - 2\Delta t Df_i^n.
\label{eqn::coarseLF}
\end{equation}
However, we would like the result from the coarse integrator to be close to that of the fine integrator. The latter is given, in our case, by applying \eqref{eqn::fineLF} twice (computing the nonlinear analogue of $\Phi^2$, as in \eqref{eqn::schurComplement} with $m=2$), which results in
\begin{equation}
\begin{array}{rl}
  u_i^{n+2} =& E(Eu_i^n - \Delta t Df_i^n) - \Delta t Df(Eu_i^n - \Delta t Dfi^n)\\
            \approx& E^2u_i^n - \Delta t\left( EDf_i^n + D f(Eu_i^n)\right)\\
            +&\Delta t^2D\left( f'(Eu_i^n)D f_i^n\right)\\
            -&\Delta t^3D\left( f''(Eu_i^n)(Df_i^n)^2\right).
\end{array}
\label{eqn::fineLFx2}
\end{equation}
Here, we Taylor-expanded the flux in the second term around $E(u_i^n)$, exploiting the fact that $\Delta t$ is a small parameter. In the case of \eqref{eqn::burgers}, this formula is exact, since the flux $f$ is a polynomial of degree $2$. By directly comparing \eqref{eqn::fineLFx2} and \eqref{eqn::coarseLF}, we see that the formula for the coarse integrator is indeed much cheaper to compute, but it also significantly differs from that of the compounded fine integrator. A way to improve the accuracy of the coarse integrator without excessively increasing its cost consists in correcting \eqref{eqn::coarseLF} so that it matches \eqref{eqn::fineLFx2} up to a certain order of $\Delta t$. This gives rise to the following integrators:
\begin{equation}
u_i^{n+1} = \left[\Psi^{(0)}(\boldsymbol{u}_i^n)\right]_i \coloneqq E^2u_i^n - 2\Delta t Df_i^n,
\label{eqn::coarseLFord0}
\end{equation}
if we aim for a zeroth-order match, leaving the rest untouched, or
\begin{equation}
u_i^{n+1} = \left[\Psi^{(1)}(\boldsymbol{u}_i^n)\right]_i \coloneqq E^2u_i^n - \Delta t\left( EDf_i^n + D f(Eu_i^n)\right),
\label{eqn::coarseLFord1}
\end{equation}
for a first-order match.
Formula \eqref{eqn::fineLFx2} refers to a coarsening factor $m=2$, but this can be easily extended to any value of $m>1$, as follows:
\begin{equation}
  u_i^{n+m} = E^{m}u_i^n - \Delta t \sum_{i=1}^{m}E^{i-1}Df(E^{m-i}u_i^n)+\mathcal{O}(\Delta t^2),
  \label{eqn::fineLFxM}
\end{equation}
and the corresponding \eqref{eqn::coarseLFord0} and \eqref{eqn::coarseLFord1} can be modified accordingly.

Unfortunately, this approach presents a number of downsides. First of all, it carries an increased computational cost with respect to rediscretisation, as the stencil of the operators involved grows linearly with $m$. For zeroth-order matching, this simply translates into a larger number of vector additions, so that the cost is still contained, provided we refrain from using aggressive coarsening strategies. For first-order matching, though, the number of flux evaluations per time-step increases as well, which makes the method far less appealing.
Secondly, this approach is limited in its application, as it requires an explicit formula for the repeated application of the fine solver \eqref{eqn::fineLFxM}: this needs to be both readily computable, and prone to simplifications without having to resort to evaluating intermediate states. All these requirements drastically limit the type of problems that can be addressed, as well as the possible schemes we can employ. In fact, we have to resort to a low-order, highly dissipative method: notice that \eqref{eqn::fineLF} corresponds to the Lax-Friedrichs flux \eqref{eqn::LFflux} with a choice of an even larger parameter $\alpha_{LF}=\Delta x/\Delta t$. For these reasons, we still recommend higher-order schemes, which in our tests behave reasonably well in most regimes. Nonetheless, the results in \secref{sec::results} show the validity of this method, and highlight the importance of accurately mimicking the action of the fine solver at the coarser levels to achieve fast convergence.

\section{Numerical results}
\label{sec::results}
In this section, we discuss how the convergence behaviour of MGRIT is impacted by factors such as the parameters in the multigrid algorithm, and the type of integrators chosen. The code used for the simulations is publicly available at \cite{myCode}.

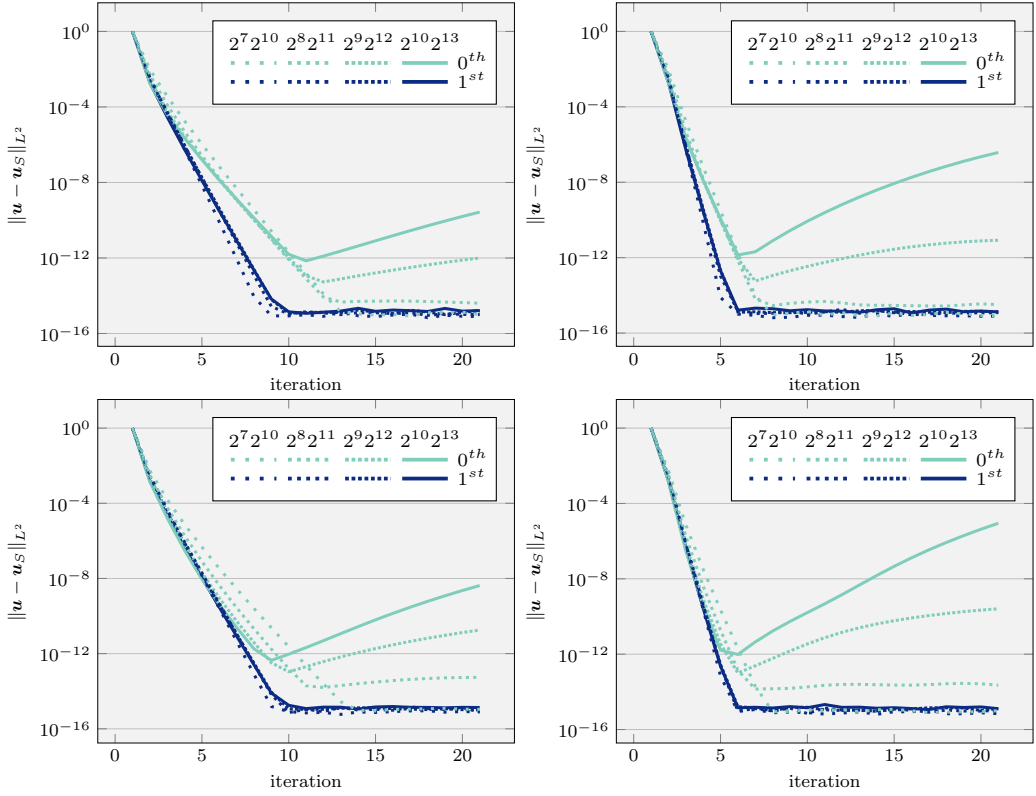
\begin{figure*}[ht!]
\centering
\ifcompileFigs
  \begin{minipage}[b]{0.49\textwidth}
    \centering
    \begin{tikzpicture}[baseline,scale=0.8]
      \begin{semilogyaxis}[tick label style={font=\small},
                   unbounded coords=jump,           
                   xlabel= iteration,
                   ylabel= $\|\boldsymbol{u}-\boldsymbol{u}_S\|_{L^2}$,
                   ymajorgrids= true,
                   axis background/.style={fill=gray!10},
                   cycle list = {{YlGnBu-L},{YlGnBu-F}} ]   
        \addplot+[                line width= 0.50mm] table[x index=0, y index=15] {data/testBurgers_ICstat_V-cycle_nlvl5_CFLmax0p95_Fsmooth.dat}; \label{fig::MGRITfineMatchSmooth::V15};
        \addplot+[                line width= 0.50mm] table[x index=0, y index=14] {data/testBurgers_ICstat_V-cycle_nlvl5_CFLmax0p95_Fsmooth.dat}; \label{fig::MGRITfineMatchSmooth::V14};
        \addplot+[densely dotted, line width= 0.50mm] table[x index=0, y index=12] {data/testBurgers_ICstat_V-cycle_nlvl5_CFLmax0p95_Fsmooth.dat}; \label{fig::MGRITfineMatchSmooth::V12};
        \addplot+[densely dotted, line width= 0.50mm] table[x index=0, y index=11] {data/testBurgers_ICstat_V-cycle_nlvl5_CFLmax0p95_Fsmooth.dat}; \label{fig::MGRITfineMatchSmooth::V11};
        \addplot+[dotted,         line width= 0.50mm] table[x index=0, y index=9 ] {data/testBurgers_ICstat_V-cycle_nlvl5_CFLmax0p95_Fsmooth.dat}; \label{fig::MGRITfineMatchSmooth::V9};
        \addplot+[dotted,         line width= 0.50mm] table[x index=0, y index=8 ] {data/testBurgers_ICstat_V-cycle_nlvl5_CFLmax0p95_Fsmooth.dat}; \label{fig::MGRITfineMatchSmooth::V8};
        \addplot+[loosely dotted, line width= 0.50mm] table[x index=0, y index=6 ] {data/testBurgers_ICstat_V-cycle_nlvl5_CFLmax0p95_Fsmooth.dat}; \label{fig::MGRITfineMatchSmooth::V6};
        \addplot+[loosely dotted, line width= 0.50mm] table[x index=0, y index=5 ] {data/testBurgers_ICstat_V-cycle_nlvl5_CFLmax0p95_Fsmooth.dat}; \label{fig::MGRITfineMatchSmooth::V5};
        \coordinate (legend) at (axis description cs:0.72,0.73);
      \end{semilogyaxis}
      \matrix [scale=0.6,
              fill=white,
              draw = black,
              matrix of nodes,
              anchor=north east,
              nodes={font=\scriptsize,text width=5.4mm,align=center,text height=1.8mm},
              ] at (legend) {    
              $2^72^{10}$                    & $2^82^{11}$                    & $2^92^{12}$                     & $2^{10}2^{13}$                 &          \\[-4pt]
              \ref{fig::MGRITfineMatchSmooth::V5}  & \ref{fig::MGRITfineMatchSmooth::V8}  & \ref{fig::MGRITfineMatchSmooth::V11}  & \ref{fig::MGRITfineMatchSmooth::V14} & |[text width= 3mm]|$0^{th}$ \\[-4pt]
              \ref{fig::MGRITfineMatchSmooth::V6}  & \ref{fig::MGRITfineMatchSmooth::V9}  & \ref{fig::MGRITfineMatchSmooth::V12}  & \ref{fig::MGRITfineMatchSmooth::V15} & |[text width= 3mm]|$1^{st}$ \\
          };
    \end{tikzpicture}
  \end{minipage}                      
  \hfill
  \begin{minipage}[b]{0.49\textwidth}
    \centering
    \begin{tikzpicture}[baseline,scale=0.8]
      \begin{semilogyaxis}[tick label style={font=\small},
                   unbounded coords=jump,           
                   xlabel= iteration,
                   ylabel= $\|\boldsymbol{u}-\boldsymbol{u}_S\|_{L^2}$,
                   ymajorgrids= true,
                   axis background/.style={fill=gray!10},   
                   cycle list = {{YlGnBu-L},{YlGnBu-F}} ]  
        \addplot+[                line width= 0.50mm] table[x index=0, y index=15] {data/testBurgers_ICstat_F-cycle_nlvl5_CFLmax0p95_Fsmooth.dat}; \label{fig::MGRITfineMatchSmooth::F15};
        \addplot+[                line width= 0.50mm] table[x index=0, y index=14] {data/testBurgers_ICstat_F-cycle_nlvl5_CFLmax0p95_Fsmooth.dat}; \label{fig::MGRITfineMatchSmooth::F14};
        \addplot+[densely dotted, line width= 0.50mm] table[x index=0, y index=12] {data/testBurgers_ICstat_F-cycle_nlvl5_CFLmax0p95_Fsmooth.dat}; \label{fig::MGRITfineMatchSmooth::F12};
        \addplot+[densely dotted, line width= 0.50mm] table[x index=0, y index=11] {data/testBurgers_ICstat_F-cycle_nlvl5_CFLmax0p95_Fsmooth.dat}; \label{fig::MGRITfineMatchSmooth::F11};
        \addplot+[dotted,         line width= 0.50mm] table[x index=0, y index=9 ] {data/testBurgers_ICstat_F-cycle_nlvl5_CFLmax0p95_Fsmooth.dat}; \label{fig::MGRITfineMatchSmooth::F9};
        \addplot+[dotted,         line width= 0.50mm] table[x index=0, y index=8 ] {data/testBurgers_ICstat_F-cycle_nlvl5_CFLmax0p95_Fsmooth.dat}; \label{fig::MGRITfineMatchSmooth::F8};
        \addplot+[loosely dotted, line width= 0.50mm] table[x index=0, y index=6 ] {data/testBurgers_ICstat_F-cycle_nlvl5_CFLmax0p95_Fsmooth.dat}; \label{fig::MGRITfineMatchSmooth::F6};
        \addplot+[loosely dotted, line width= 0.50mm] table[x index=0, y index=5 ] {data/testBurgers_ICstat_F-cycle_nlvl5_CFLmax0p95_Fsmooth.dat}; \label{fig::MGRITfineMatchSmooth::F5};
        \coordinate (legend) at (axis description cs:0.72,0.73);
      \end{semilogyaxis}
      \matrix [scale=0.6,
              fill=white,
              draw = black,
              matrix of nodes,
              anchor=north east,
              nodes={font=\scriptsize,text width=5.4mm,align=center,text height=1.8mm},
              ] at (legend) {            
              $2^72^{10}$                    & $2^82^{11}$                    & $2^92^{12}$                     & $2^{10}2^{13}$                 &          \\[-4pt]
              \ref{fig::MGRITfineMatchSmooth::F5}  & \ref{fig::MGRITfineMatchSmooth::F8}  & \ref{fig::MGRITfineMatchSmooth::F11}  & \ref{fig::MGRITfineMatchSmooth::F14} & |[text width= 3mm]|$0^{th}$ \\[-4pt]
              \ref{fig::MGRITfineMatchSmooth::F6}  & \ref{fig::MGRITfineMatchSmooth::F9}  & \ref{fig::MGRITfineMatchSmooth::F12}  & \ref{fig::MGRITfineMatchSmooth::F15} & |[text width= 3mm]|$1^{st}$ \\
          };
    \end{tikzpicture}
  \end{minipage}
  \par
  \begin{minipage}[b]{0.49\textwidth}
    \centering
    \begin{tikzpicture}[baseline,scale=0.8]
      \begin{semilogyaxis}[tick label style={font=\small},
                   unbounded coords=jump,           
                   xlabel= iteration,
                   ylabel= $\|\boldsymbol{u}-\boldsymbol{u}_S\|_{L^2}$,
                   ymajorgrids= true,
                   axis background/.style={fill=gray!10},
                   cycle list = {{YlGnBu-L},{YlGnBu-F}} ]   
        \addplot+[                line width= 0.50mm] table[x index=0, y index=15] {data/testBurgers_ICmov_V-cycle_nlvl5_CFLmax0p95_Fsmooth.dat}; \label{fig::MGRITfineMatch::V15};
        \addplot+[                line width= 0.50mm] table[x index=0, y index=14] {data/testBurgers_ICmov_V-cycle_nlvl5_CFLmax0p95_Fsmooth.dat}; \label{fig::MGRITfineMatch::V14};
        \addplot+[densely dotted, line width= 0.50mm] table[x index=0, y index=12] {data/testBurgers_ICmov_V-cycle_nlvl5_CFLmax0p95_Fsmooth.dat}; \label{fig::MGRITfineMatch::V12};
        \addplot+[densely dotted, line width= 0.50mm] table[x index=0, y index=11] {data/testBurgers_ICmov_V-cycle_nlvl5_CFLmax0p95_Fsmooth.dat}; \label{fig::MGRITfineMatch::V11};
        \addplot+[dotted,         line width= 0.50mm] table[x index=0, y index=9 ] {data/testBurgers_ICmov_V-cycle_nlvl5_CFLmax0p95_Fsmooth.dat}; \label{fig::MGRITfineMatch::V9};
        \addplot+[dotted,         line width= 0.50mm] table[x index=0, y index=8 ] {data/testBurgers_ICmov_V-cycle_nlvl5_CFLmax0p95_Fsmooth.dat}; \label{fig::MGRITfineMatch::V8};
        \addplot+[loosely dotted, line width= 0.50mm] table[x index=0, y index=6 ] {data/testBurgers_ICmov_V-cycle_nlvl5_CFLmax0p95_Fsmooth.dat}; \label{fig::MGRITfineMatch::V6};
        \addplot+[loosely dotted, line width= 0.50mm] table[x index=0, y index=5 ] {data/testBurgers_ICmov_V-cycle_nlvl5_CFLmax0p95_Fsmooth.dat}; \label{fig::MGRITfineMatch::V5};
        \coordinate (legend) at (axis description cs:0.72,0.73);
      \end{semilogyaxis}
      \matrix [scale=0.6,
              fill=white,
              draw = black,
              matrix of nodes,
              anchor=north east,
              nodes={font=\scriptsize,text width=5.4mm,align=center,text height=1.8mm},
              ] at (legend) {
              $2^72^{10}$                    & $2^82^{11}$                    & $2^92^{12}$                     & $2^{10}2^{13}$                 &          \\[-4pt]
              \ref{fig::MGRITfineMatch::V5}  & \ref{fig::MGRITfineMatch::V8}  & \ref{fig::MGRITfineMatch::V11}  & \ref{fig::MGRITfineMatch::V14} & |[text width= 3mm]|$0^{th}$ \\[-4pt]
              \ref{fig::MGRITfineMatch::V6}  & \ref{fig::MGRITfineMatch::V9}  & \ref{fig::MGRITfineMatch::V12}  & \ref{fig::MGRITfineMatch::V15} & |[text width= 3mm]|$1^{st}$ \\
          };
    \end{tikzpicture}
  \end{minipage}                      
  \hfill
  \begin{minipage}[b]{0.49\textwidth}
    \centering
    \begin{tikzpicture}[baseline,scale=0.8]
      \begin{semilogyaxis}[tick label style={font=\small},
                   unbounded coords=jump,           
                   xlabel= iteration,
                   ylabel= $\|\boldsymbol{u}-\boldsymbol{u}_S\|_{L^2}$,
                   ymajorgrids= true,
                   axis background/.style={fill=gray!10},
                   cycle list = {{YlGnBu-L},{YlGnBu-F}} ]  
        \addplot+[                line width= 0.50mm] table[x index=0, y index=15] {data/testBurgers_ICmov_F-cycle_nlvl5_CFLmax0p95_Fsmooth.dat}; \label{fig::MGRITfineMatch::F15};
        \addplot+[                line width= 0.50mm] table[x index=0, y index=14] {data/testBurgers_ICmov_F-cycle_nlvl5_CFLmax0p95_Fsmooth.dat}; \label{fig::MGRITfineMatch::F14};
        \addplot+[densely dotted, line width= 0.50mm] table[x index=0, y index=12] {data/testBurgers_ICmov_F-cycle_nlvl5_CFLmax0p95_Fsmooth.dat}; \label{fig::MGRITfineMatch::F12};
        \addplot+[densely dotted, line width= 0.50mm] table[x index=0, y index=11] {data/testBurgers_ICmov_F-cycle_nlvl5_CFLmax0p95_Fsmooth.dat}; \label{fig::MGRITfineMatch::F11};
        \addplot+[dotted,         line width= 0.50mm] table[x index=0, y index=9 ] {data/testBurgers_ICmov_F-cycle_nlvl5_CFLmax0p95_Fsmooth.dat}; \label{fig::MGRITfineMatch::F9};
        \addplot+[dotted,         line width= 0.50mm] table[x index=0, y index=8 ] {data/testBurgers_ICmov_F-cycle_nlvl5_CFLmax0p95_Fsmooth.dat}; \label{fig::MGRITfineMatch::F8};
        \addplot+[loosely dotted, line width= 0.50mm] table[x index=0, y index=6 ] {data/testBurgers_ICmov_F-cycle_nlvl5_CFLmax0p95_Fsmooth.dat}; \label{fig::MGRITfineMatch::F6};
        \addplot+[loosely dotted, line width= 0.50mm] table[x index=0, y index=5 ] {data/testBurgers_ICmov_F-cycle_nlvl5_CFLmax0p95_Fsmooth.dat}; \label{fig::MGRITfineMatch::F5};
        \coordinate (legend) at (axis description cs:0.72,0.73);
      \end{semilogyaxis}
      \matrix [scale=0.6,
              fill=white,
              draw = black,
              matrix of nodes,
              anchor=north east,
              nodes={font=\scriptsize,text width=5.4mm,align=center,text height=1.8mm},
              ] at (legend) {
              $2^72^{10}$                    & $2^82^{11}$                    & $2^92^{12}$                     & $2^{10}2^{13}$                 &          \\[-4pt]
              \ref{fig::MGRITfineMatch::F5}  & \ref{fig::MGRITfineMatch::F8}  & \ref{fig::MGRITfineMatch::F11}  & \ref{fig::MGRITfineMatch::F14} & |[text width= 3mm]|$0^{th}$ \\[-4pt]
              \ref{fig::MGRITfineMatch::F6}  & \ref{fig::MGRITfineMatch::F9}  & \ref{fig::MGRITfineMatch::F12}  & \ref{fig::MGRITfineMatch::F15} & |[text width= 3mm]|$1^{st}$ \\
          };
    \end{tikzpicture}
  \end{minipage}
\else
  \begin{minipage}[b]{0.49\textwidth}
    \centering
    \includegraphics{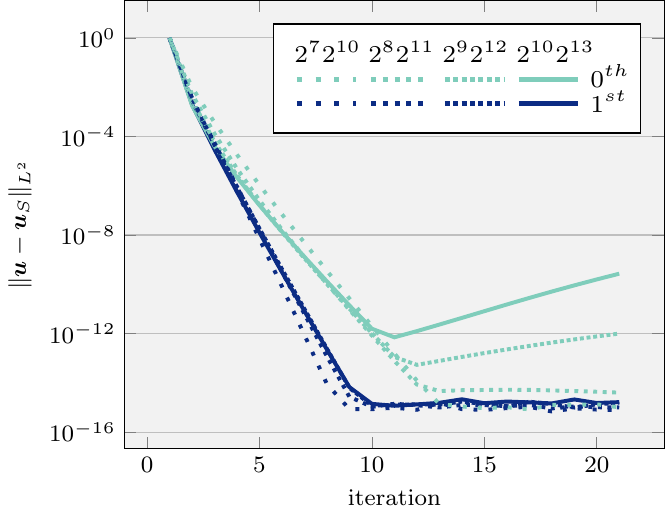}
  \end{minipage}                      
  \hfill
  \begin{minipage}[b]{0.49\textwidth}
    \centering
    \includegraphics{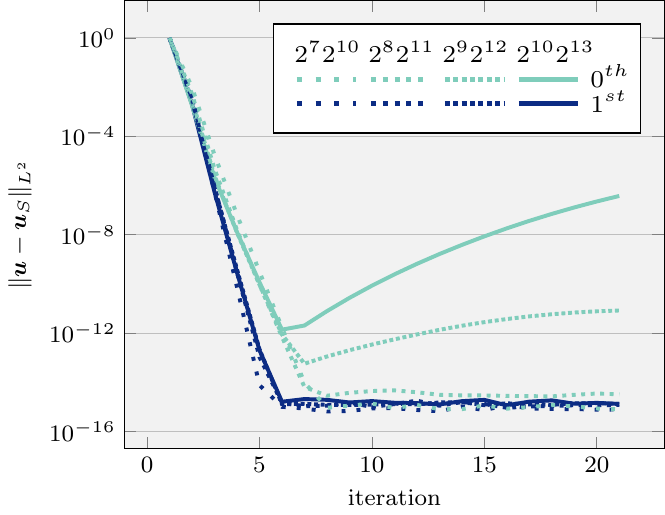}
  \end{minipage}
  \par
  \begin{minipage}[b]{0.49\textwidth}
    \centering
    \includegraphics{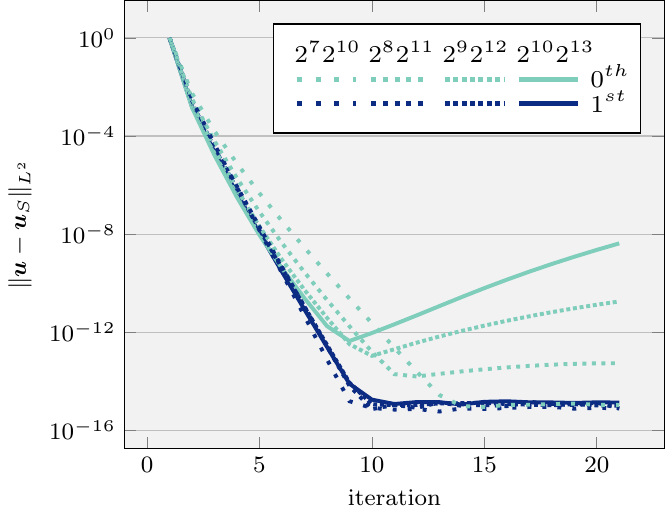}
  \end{minipage}                      
  \hfill
  \begin{minipage}[b]{0.49\textwidth}
    \centering
    \includegraphics{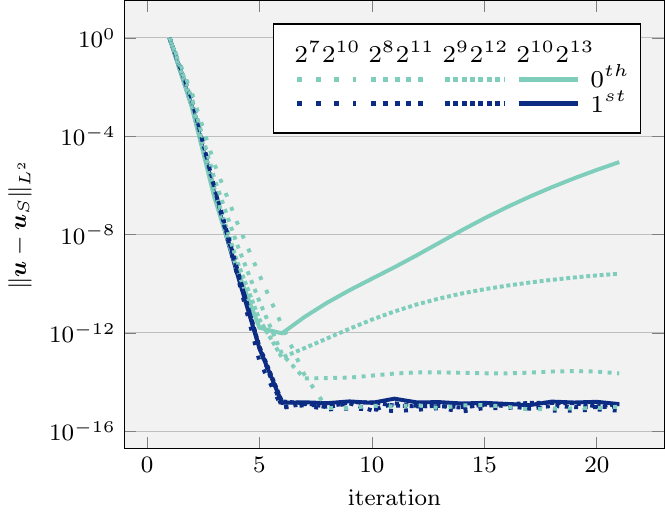}
  \end{minipage}
\fi
\caption[Error convergence for fine integrator matching]{Impact of grid refinement and degree of fine-integrator matching on the evolution of the error between the MGRIT iterates and the serial solution $\boldsymbol{u}_S$. The relative $L^2$ norm in space-time is considered. MGRIT parameters: F-smoothing, $5$ levels, coarsening factor $m=2$, V-cycle (left), and F-cycle (right). Domain: $L=1$, $T=0.475$. More densely dotted graphs correspond to finer meshes: the numbers of nodes $N_xN_t$ at the finest level are reported in the legend. The time-steppers introduced in \secref{sec::fineMatching} are employed: FE with Lax-Friedrichs flux \eqref{eqn::fineLF} for the finest level, and the corresponding zeroth- \eqref{eqn::coarseLFord0} and first-order matching \eqref{eqn::coarseLFord1}, with varying coarsening factors, for the coarse levels (green and blue lines, respectively). Initial conditions: stationary shock (top, $u(x,0)=\bar{u}^0_S$), and moving shock (bottom, $u(x,0)=\bar{u}^0_M$). Results from pure rediscretisation are not reported, as the error diverges after the very first iteration. Data generated with \texttt{testMGRIT\_Burgers\_Matching.m} \cite{myCode}}
\label{fig::MGRITfineMatch}
\end{figure*}

\paragraph{Fine integrator matching} As a first test, we check the effectiveness of the time-stepper proposed in \secref{sec::fineMatching} applied to Burgers' equation: the convergence results for a number of runs are reported in \figref{fig::MGRITfineMatch}. We choose two different initial conditions: first, the simple sinusoidal wave $\bar{u}^0_S = \sin(2\pi x/L)$, for which the solution develops a \emph{stationary shock} positioned at $x_S=L/2$, at the \emph{breaking time} $t_S=-1/\min_x(\bar{u}^0_S(x))'=L/(2\pi)$; second, its translation $\bar{u}^0_M = ( 1+\sin(2\pi x/L) )/2$, which instead develops a moving shock. These two choices allow us to investigate the impact of having to track a discontinuity which is not aligned with the grid.

\begin{figure}[hb!]
\centering
\ifcompileFigs
  \begin{tikzpicture}[baseline,scale=0.8]
    \begin{semilogyaxis}[tick label style={font=\small},
                 unbounded coords=jump,           
                 xlabel= iteration,
                 ylabel= $\|\boldsymbol{u}-\boldsymbol{u}_S\|_{L^2}$,
                 ymajorgrids= true,
                 axis background/.style={fill=gray!10},
                 cycle list = {{YlGnBu-L},{YlGnBu-F}} ]  
      \addplot+[         line width= 0.50mm] table[x index=0, y index=5 ] {data/testBurgers_ICstat_V-cycle_nlvl5_CFLmax0p95_Fsmooth.dat}; \label{fig::MGRITrestrictionGuess::0};
      \addplot+[         line width= 0.50mm] table[x index=0, y index=5 ] {data/testBurgers_ICstat_V-cycle_nlvl5_CFLmax0p95_Fsmooth_BadRestrictionGuess.dat}; \label{fig::MGRITrestrictionGuess::1};
      \addplot+[ dashed, line width= 0.50mm] table[x index=0, y index=5 ] {data/testBurgers_ICstat_V-cycle_nlvl5_CFLmax0p95_F-CFsmooth.dat}; \label{fig::MGRITrestrictionGuess::2};
      \addplot+[ dashed, line width= 0.50mm] table[x index=0, y index=5 ] {data/testBurgers_ICstat_V-cycle_nlvl5_CFLmax0p95_F-CFsmooth_BadRestrictionGuess.dat}; \label{fig::MGRITrestrictionGuess::3};
      \coordinate (legend) at (axis description cs:0.725,0.735);
    \end{semilogyaxis}
    \matrix [scale=0.6,
            fill=white,
            draw = black,
            matrix of nodes,
            anchor=north east,
            nodes={font=\scriptsize,text width=6.4mm,align=center,text height=1.4mm},
            ] at (legend) {
            $\text{F}$                          & $\text{FCF}$                        &                                                                 \\[-4pt]
            \ref{fig::MGRITrestrictionGuess::1} & \ref{fig::MGRITrestrictionGuess::3} & |[text width= 17mm]|$\text{Injection}$                          \\[-4pt]
            \ref{fig::MGRITrestrictionGuess::0} & \ref{fig::MGRITrestrictionGuess::2} & |[text width= 17mm]|$\text{Formula \eqref{eqn::restriction::u}}$\\[-4pt]
        };
  \end{tikzpicture}
\else
  \includegraphics{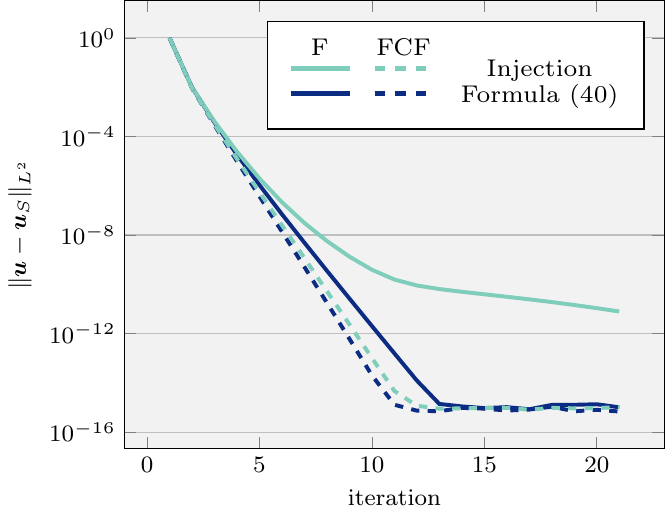}
\fi
\caption[Error convergence for different initial guesses]{Effect of modifying the initial guess for the coarse solution on the error convergence of MGRIT: comparison of error convergence between classical injection (green) and formula \eqref{eqn::restriction::u} (blue). Relaxation used: F-relaxation (full lines) and FCF-relaxation (dashed lines); otherwise, same parameters as in the top-left graph in \figref{fig::MGRITfineMatch}. Only results for $N_xN_t=2^{7}2^{10}$ (zeroth-order) are reported, but finer grids give comparable results}
\label{fig::MGRITrestrictionGuess}
\end{figure}
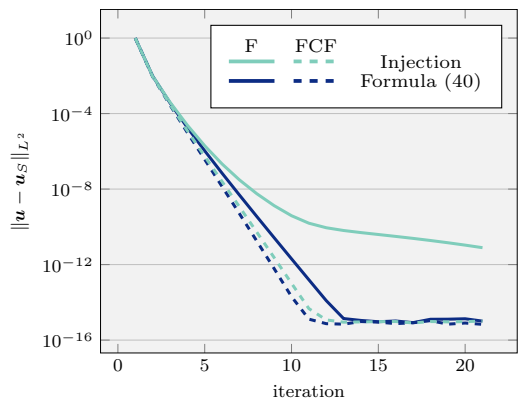

First of all, we observe that the performance definitely improves if the coarse solver is chosen to match the fine solver up to a higher order: in fact, the blue curves in \figref{fig::MGRITfineMatch} (first-order matching, \eqref{eqn::coarseLFord1}) lie consistently below the green ones (zeroth-order matching, \eqref{eqn::coarseLFord0}). This is in line with expectations, since we are employing a more accurate solver. Still, the performance of zeroth-order matching remains competitive, achieving very effective reductions in error over the first iterations. Somewhat surprisingly, the zeroth-order matching gets slightly better initial error reductions in the moving-shock scenario than for the stationary shock. Unfortunately, though, this scheme clearly fails to effectively damp some components of the error, as we can see from the fact that the corresponding error starts rising again after 5-10 iterations on finer meshes. This negative effect seems to be heightened as we refine the grid further, with the error starting to increase earlier. Nonetheless, in all the tests considered this approach behaves much better than na\"ive rediscretisation, for which the error blows up already after the very first iteration and is, hence, not reported in the graph. Moreover, the convergence behaviour remains unchanged as we further refine our grids: the different lines superimpose almost perfectly. Modifying the type of multigrid cycle does not vary the nature of these observations, but as expected an F-cycle shows steeper convergence plots than a V-cycle. We have found that modifying the type of relaxation from an F- to an FCF-smoother, instead, does not improve convergence dramatically, so long as the initial guess for the solution at the coarse nodes is picked as in \eqref{eqn::restriction::u}. An example of this is shown in \figref{fig::MGRITrestrictionGuess}, where a comparison of the error evolution of MGRIT is shown when using simple injection as a restriction operator versus formula \eqref{eqn::restriction::u}. The latter consistently shows better (or at worst comparable) results, particularly if MGRIT is equipped with simple multigrid components, such as V-cycle and F-relaxation as opposed to F-cycle and FCF-relaxation. 

\begin{figure*}[ht!]
\centering
\ifcompileFigs
  \begin{minipage}[b]{0.49\textwidth}
    \centering
    \begin{tikzpicture}[baseline,scale=0.8]
      \begin{semilogyaxis}[tick label style={font=\small},
                   xlabel= iteration,
                   ylabel= $\|\boldsymbol{u}-\boldsymbol{u}_S\|_{L^2}$,
                   ymajorgrids= true,
                   axis background/.style={fill=gray!10},
                   legend pos=north east,
                   cycle list = {{YlGnBu-D}, {YlGnBu-F}, {YlGnBu-H}, {YlGnBu-L}} ]
        \addplot+[        line width=0.5mm] table[x index=0, y index=6] {data/testBurgers_ICstat_V-cycle_nlvl5_CFLmax7p6_Fsmooth.dat}; 
        \addplot+[        line width=0.5mm] table[x index=0, y index=6] {data/testBurgers_ICstat_V-cycle_nlvl5_CFLmax3p8_Fsmooth.dat}; 
        \addplot+[        line width=0.5mm] table[x index=0, y index=6] {data/testBurgers_ICstat_V-cycle_nlvl5_CFLmax1p9_Fsmooth.dat}; 
        \addplot+[dashed, line width=0.5mm] table[x index=0, y index=6] {data/testBurgers_ICstat_V-cycle_nlvl5_CFLmax0p95_Fsmooth.dat};
        \addlegendentry{$c=7.6$}
        \addlegendentry{$c=3.8$}
        \addlegendentry{$c=1.9$}
        \addlegendentry{$\boldsymbol{c=.95}$}
        \addlegendentry{$c=.47$}
        \addlegendentry{$c=.24$}
      \end{semilogyaxis}
    \end{tikzpicture}
  \end{minipage}
  \hfill
  \begin{minipage}[b]{0.49\textwidth}
    \centering
    \begin{tikzpicture}[baseline,scale=0.8]
      \begin{semilogyaxis}[tick label style={font=\small},
                   unbounded coords=jump,           
                   xlabel= iteration,
                   ylabel= $\|\boldsymbol{u}-\boldsymbol{u}_S\|_{L^2}$,
                   ymajorgrids= true,
                   ymax = 1000,
                   axis background/.style={fill=gray!10},
                   legend pos=north east,
                   cycle list = {{YlGnBu-D}, {YlGnBu-F}, {YlGnBu-H}, {YlGnBu-L}} ]
        \addplot+[        line width=0.50mm] table[x index=0, y index=1 ] {data/MGRIT_CFL_Burgers1D_WENOROE_Nlvl5_NX128_maxNT4096_maxS7.dat};
        \addplot+[dashed, line width=0.50mm] table[x index=0, y index=5 ] {data/MGRIT_CFL_Burgers1D_WENOROE_Nlvl5_NX128_maxNT4096_maxS7.dat};
        \addplot+[        line width=0.50mm] table[x index=0, y index=9 ] {data/MGRIT_CFL_Burgers1D_WENOROE_Nlvl5_NX128_maxNT4096_maxS7.dat};
        \addplot+[        line width=0.50mm] table[x index=0, y index=13] {data/MGRIT_CFL_Burgers1D_WENOROE_Nlvl5_NX128_maxNT4096_maxS7.dat};
        \addlegendentry{$c=1.9$}
        \addlegendentry{$\boldsymbol{c=.95}$}
        \addlegendentry{$c=.48$}
        \addlegendentry{$c=.24$}
      \end{semilogyaxis}
    \end{tikzpicture}
  \end{minipage}
  \par
  \begin{minipage}[b]{0.49\textwidth}
    \centering
    \begin{tikzpicture}[baseline,scale=0.8]
      \begin{semilogyaxis}[tick label style={font=\small},
                   unbounded coords=jump,           
                   xlabel= iteration,
                   ylabel= $\|\boldsymbol{u}-\boldsymbol{u}_S\|_{L^2}$,
                   ymajorgrids= true,
                   ymax = 1000,
                   axis background/.style={fill=gray!10},
                   legend pos=north east,
                   cycle list = {{YlGnBu-F}, {YlGnBu-H}, {YlGnBu-L}} ]
        \addplot+[dashed, line width=0.50mm] table[x index=0, y index=7 ] {data/MGRIT_CFL_Burgers1D_WENOROE_Nlvl5_NX128_maxNT4096_maxS7_FE.dat};
        \addplot+[        line width=0.50mm] table[x index=0, y index=11] {data/MGRIT_CFL_Burgers1D_WENOROE_Nlvl5_NX128_maxNT4096_maxS7_FE.dat};
        \addplot+[        line width=0.50mm] table[x index=0, y index=15] {data/MGRIT_CFL_Burgers1D_WENOROE_Nlvl5_NX128_maxNT4096_maxS7_FE.dat};
        \addlegendentry{$\boldsymbol{c=.95}$}
        \addlegendentry{$c=.48$}
        \addlegendentry{$c=.24$}
      \end{semilogyaxis}
    \end{tikzpicture}
  \end{minipage}
  \hfill
  \begin{minipage}[b]{0.49\textwidth}
    \centering
    \begin{tikzpicture}[baseline,scale=0.8]
      \begin{semilogyaxis}[tick label style={font=\small},
                   unbounded coords=jump,           
                   xlabel= iteration,
                   ylabel= $\|\boldsymbol{u}-\boldsymbol{u}_S\|_{L^2}$,
                   ymajorgrids= true,
                   ymax = 1000,
                   axis background/.style={fill=gray!10},
                   legend pos=north east,
                   cycle list = {{YlGnBu-D}, {YlGnBu-F}, {YlGnBu-H}, {YlGnBu-L}} ]
        \addplot+[        line width=0.50mm] table[x index=0, y index=3 ] {data/MGRIT_CFL_Burgers1D_WENOROE_Nlvl5_NX128_maxNT4096_maxS7.dat};
        \addplot+[dashed, line width=0.50mm] table[x index=0, y index=7 ] {data/MGRIT_CFL_Burgers1D_WENOROE_Nlvl5_NX128_maxNT4096_maxS7.dat};
        \addplot+[        line width=0.50mm] table[x index=0, y index=11] {data/MGRIT_CFL_Burgers1D_WENOROE_Nlvl5_NX128_maxNT4096_maxS7.dat};
        \addplot+[        line width=0.50mm] table[x index=0, y index=15] {data/MGRIT_CFL_Burgers1D_WENOROE_Nlvl5_NX128_maxNT4096_maxS7.dat};
        \addlegendentry{$c=1.9$}
        \addlegendentry{$\boldsymbol{c=.95}$}
        \addlegendentry{$c=.48$}
        \addlegendentry{$c=.24$}
      \end{semilogyaxis}
    \end{tikzpicture}
  \end{minipage}
\else
\begin{minipage}[b]{0.49\textwidth}
    \centering
    \includegraphics{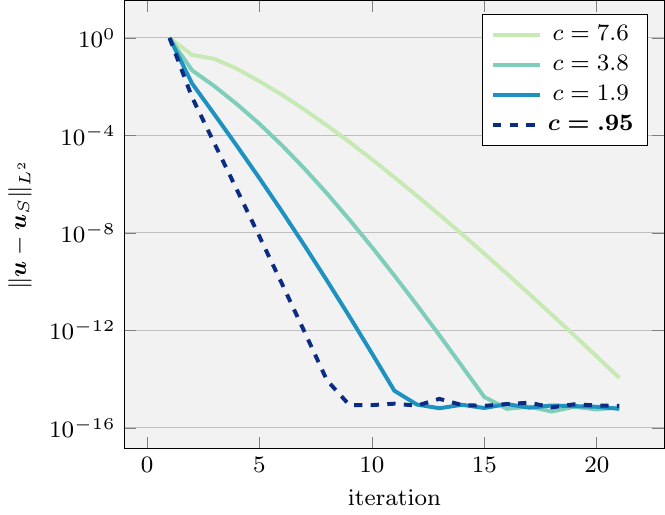}
  \end{minipage}
  \hfill
  \begin{minipage}[b]{0.49\textwidth}
    \centering
    \includegraphics{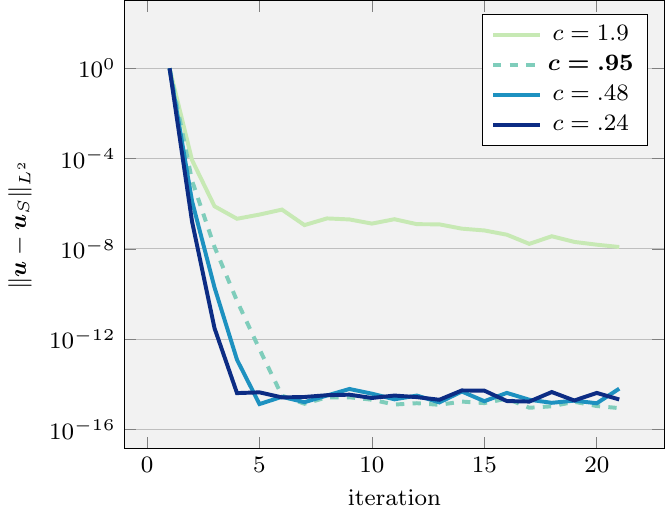}
  \end{minipage}
  \par
  \begin{minipage}[b]{0.49\textwidth}
    \centering
    \includegraphics{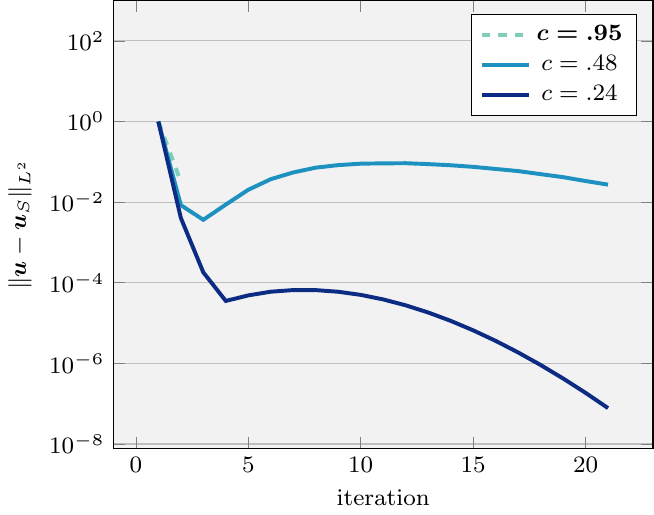}
  \end{minipage}
  \hfill
  \begin{minipage}[b]{0.49\textwidth}
    \centering
    \includegraphics{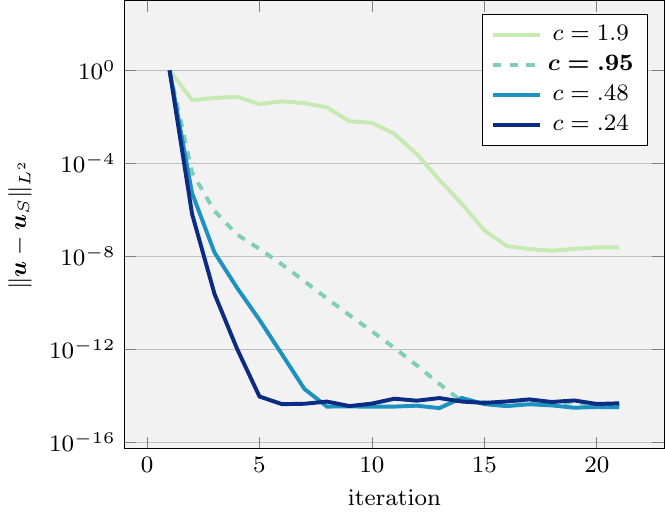}
  \end{minipage}
\fi
\caption[Error convergence for different CFL numbers]{Influence of increasing the CFL number $c$ at the coarsest level on the convergence behaviour of MGRIT applied to Burgers' equation. $N_x=2^7$, while $N_t$ is progressively halved to increase $c$. The other parameters are the same as in the top-left graph in \figref{fig::MGRITfineMatch}, but different schemes are used in each plot. Top-left: first-order matching and FE; top-right: WENO reconstruction of order $s=1$ and SSPRK3; bottom-left: WENO ($s=5$) and FE; bottom-right: WENO ($s=5$) and SSPRK3; the last three use a Roe flux and rediscretisation at all levels. Missing values are diverging (notice in particular the truncated $c=0.95$ plot in the bottom-left graph). Data generated using \texttt{testMGRIT\_highOrderCFL.m} \cite{myCode}}
\label{fig::MGRITcfl}
\end{figure*}
Rather, the factor playing the main role in determining the performance of the algorithm lies in how faithfully the \emph{Courant-Friedrichs-Lewy} CFL condition \cite[Chap.~10.6]{leveque} is respected. 
In \figref{fig::MGRITfineMatch}, the number of nodes $N_t$ and $N_x$ are chosen to guarantee that the CFL number on the \emph{coarsest} temporal mesh is $\sim0.95$ in all cases. This also implies that the solution is \emph{de facto} over-resolved, since usually one aims for a CFL number of $\sim1$ at the \emph{finest} level. The impact of relaxing the coarse-grid CFL condition can be clearly seen in the top-left graph of \figref{fig::MGRITcfl}: convergence degrades noticeably, even using high-order matching.

\paragraph{High-order schemes}
While the integrator proposed in \secref{sec::fineMatching} seems effective, it remains very dissipative. On the one hand, it is well-documented (see for example \cite{MGRITconv}) that smoothing effects improve convergence of MGRIT; on the other, though, artificial diffusion is very undesirable in the model problems considered, as it results in a heavy smearing of the shocks, which are the characteristic feature of the solutions of conservation laws. The choice of high-order space reconstructions is instead preferable, as they offer the possibility to preserve such discontinuities sharply, and an equally high-order time-discretisation is consequently requested. Investigating the behaviour of the MGRIT algorithm when used in conjunction with these schemes is hence of relevance, as it would more faithfully represent the setup of a real-world application.
For this reason, we proceed to applying MGRIT to the test cases introduced in \secref{sec::modelProblems}, for discretisations using the Lax-Friedrichs and Roe definitions of numerical fluxes, \eqref{eqn::LFflux} and \eqref{eqn::roeflux}, and employing WENO reconstructions at various orders of accuracy $s$.
\begin{figure*}[ht!]
\centering
\ifcompileFigs
  \begin{minipage}[b]{0.49\textwidth}
    \centering
    \begin{tikzpicture}[baseline,scale=0.8]
      \begin{semilogyaxis}[tick label style={font=\small},
                 xlabel= iteration,
                 ylabel= $\|\boldsymbol{u}-\boldsymbol{u}_S\|_{L^2}$,
                 ymajorgrids= true,
                 axis background/.style={fill=gray!10},
                 legend pos=north east,
                 cycle list = {{YlGnBu-L}, {YlGnBu-H}, {YlGnBu-F}, {YlGnBu-D}} ]
        \addplot+[line width=0.50mm] table[x index=0, y index=1] {data/MGRIT_Burgers1DWENO_LF.dat};
        \addplot+[line width=0.50mm] table[x index=0, y index=2] {data/MGRIT_Burgers1DWENO_LF.dat};
        \addplot+[line width=0.50mm] table[x index=0, y index=3] {data/MGRIT_Burgers1DWENO_LF.dat};
        \addplot+[line width=0.50mm] table[x index=0, y index=4] {data/MGRIT_Burgers1DWENO_LF.dat};
        \addlegendentry{$s=1$}
        \addlegendentry{$s=3$}
        \addlegendentry{$s=5$}
        \addlegendentry{$s=7$}
      \end{semilogyaxis}
    \end{tikzpicture}
  \end{minipage}                      
  \hfill
  \begin{minipage}[b]{0.49\textwidth}
    \centering
    \begin{tikzpicture}[baseline,scale=0.8]
      \begin{semilogyaxis}[tick label style={font=\small},
                 xlabel= iteration,
                 ylabel= $\|\boldsymbol{u}-\boldsymbol{u}_S\|_{L^2}$,
                 ymajorgrids= true,
                 axis background/.style={fill=gray!10},
                 legend pos=north east,
                 cycle list = {{YlGnBu-L}, {YlGnBu-H}, {YlGnBu-F}, {YlGnBu-D}} ]
        \addplot+[line width=0.50mm] table[x index=0, y index=1] {data/MGRIT_Burgers1DWENO_ROE.dat};
        \addplot+[line width=0.50mm] table[x index=0, y index=2] {data/MGRIT_Burgers1DWENO_ROE.dat};
        \addplot+[line width=0.50mm] table[x index=0, y index=3] {data/MGRIT_Burgers1DWENO_ROE.dat};
        \addplot+[line width=0.50mm] table[x index=0, y index=4] {data/MGRIT_Burgers1DWENO_ROE.dat};
        \addlegendentry{$s=1$}
        \addlegendentry{$s=3$}
        \addlegendentry{$s=5$}
        \addlegendentry{$s=7$}
      \end{semilogyaxis}
    \end{tikzpicture}
  \end{minipage}
  \par
  \begin{minipage}[b]{0.49\textwidth}
    \centering
    \begin{tikzpicture}[baseline,scale=0.8]
      \begin{semilogyaxis}[tick label style={font=\small},
                 xlabel= iteration,
                 ylabel= $\|\boldsymbol{u}-\boldsymbol{u}_S\|_{L^2}$,
                 ymajorgrids= true,
                 axis background/.style={fill=gray!10},
                 legend pos=north east,
                 cycle list = {{YlGnBu-L}, {YlGnBu-H}, {YlGnBu-F}, {YlGnBu-D}} ]
        \addplot+[line width=0.50mm] table[x index=0, y index=1] {data/MGRIT_shallowWater1DWENO_LF.dat};
        \addplot+[line width=0.50mm] table[x index=0, y index=2] {data/MGRIT_shallowWater1DWENO_LF.dat};
        \addplot+[line width=0.50mm] table[x index=0, y index=3] {data/MGRIT_shallowWater1DWENO_LF.dat};
        \addplot+[line width=0.50mm] table[x index=0, y index=4] {data/MGRIT_shallowWater1DWENO_LF.dat};
        \addlegendentry{$s=1$}
        \addlegendentry{$s=3$}
        \addlegendentry{$s=5$}
        \addlegendentry{$s=7$}
      \end{semilogyaxis}
    \end{tikzpicture}
  \end{minipage}                      
  \hfill
  \begin{minipage}[b]{0.49\textwidth}
    \centering
    \begin{tikzpicture}[baseline,scale=0.8]
      \begin{semilogyaxis}[tick label style={font=\small},
                 xlabel= iteration,
                 ylabel= $\|\boldsymbol{u}-\boldsymbol{u}_S\|_{L^2}$,
                 ymajorgrids= true,
                 axis background/.style={fill=gray!10},
                 legend pos=north east,
                 cycle list = {{YlGnBu-L}, {YlGnBu-H}, {YlGnBu-F}, {YlGnBu-D}} ]
        \addplot+[line width=0.50mm] table[x index=0, y index=1] {data/MGRIT_shallowWater1DWENO_ROE.dat};
        \addplot+[line width=0.50mm] table[x index=0, y index=2] {data/MGRIT_shallowWater1DWENO_ROE.dat};
        \addplot+[line width=0.50mm] table[x index=0, y index=3] {data/MGRIT_shallowWater1DWENO_ROE.dat};
        \addplot+[line width=0.50mm] table[x index=0, y index=4] {data/MGRIT_shallowWater1DWENO_ROE.dat};
        \addlegendentry{$s=1$}
        \addlegendentry{$s=3$}
        \addlegendentry{$s=5$}
        \addlegendentry{$s=7$}
      \end{semilogyaxis}
    \end{tikzpicture}
  \end{minipage}
  \par
  \begin{minipage}[b]{0.49\textwidth}
    \centering
    \begin{tikzpicture}[baseline,scale=0.8]
      \begin{semilogyaxis}[tick label style={font=\small},
                 xlabel= iteration,
                 ylabel= $\|\boldsymbol{u}-\boldsymbol{u}_S\|_{L^2}$,
                 ymajorgrids= true,
                 axis background/.style={fill=gray!10},
                 legend pos=north east,
                 cycle list = {{YlGnBu-L}, {YlGnBu-H}, {YlGnBu-F}, {YlGnBu-D}} ]
        \addplot+[line width=0.50mm] table[x index=0, y index=1] {data/MGRIT_euler1DWENO_LF.dat};
        \addplot+[line width=0.50mm] table[x index=0, y index=2] {data/MGRIT_euler1DWENO_LF.dat};
        \addplot+[line width=0.50mm] table[x index=0, y index=3] {data/MGRIT_euler1DWENO_LF.dat};
        \addplot+[line width=0.50mm] table[x index=0, y index=4] {data/MGRIT_euler1DWENO_LF.dat};
        \addlegendentry{$s=1$}
        \addlegendentry{$s=3$}
        \addlegendentry{$s=5$}
        \addlegendentry{$s=7$}
      \end{semilogyaxis}
    \end{tikzpicture}
  \end{minipage}                      
  \hfill
  \begin{minipage}[b]{0.49\textwidth}
    \centering
    \begin{tikzpicture}[baseline,scale=0.8]
      \begin{semilogyaxis}[tick label style={font=\small},
                 xlabel= iteration,
                 ylabel= $\|\boldsymbol{u}-\boldsymbol{u}_S\|_{L^2}$,
                 ymajorgrids= true,
                 axis background/.style={fill=gray!10},
                 legend pos=north east,
                 cycle list = {{YlGnBu-L}, {YlGnBu-H}, {YlGnBu-F}, {YlGnBu-D}} ]
        \addplot+[line width=0.50mm] table[x index=0, y index=1] {data/MGRIT_euler1DWENO_ROE.dat};
        \addplot+[line width=0.50mm] table[x index=0, y index=2] {data/MGRIT_euler1DWENO_ROE.dat};
        \addplot+[line width=0.50mm] table[x index=0, y index=3] {data/MGRIT_euler1DWENO_ROE.dat};
        \addplot+[line width=0.50mm] table[x index=0, y index=4] {data/MGRIT_euler1DWENO_ROE.dat};
        \addlegendentry{$s=1$}
        \addlegendentry{$s=3$}
        \addlegendentry{$s=5$}
        \addlegendentry{$s=7$}
      \end{semilogyaxis}
    \end{tikzpicture}
  \end{minipage}
\else
  \begin{minipage}[b]{0.49\textwidth}
    \centering
    \includegraphics{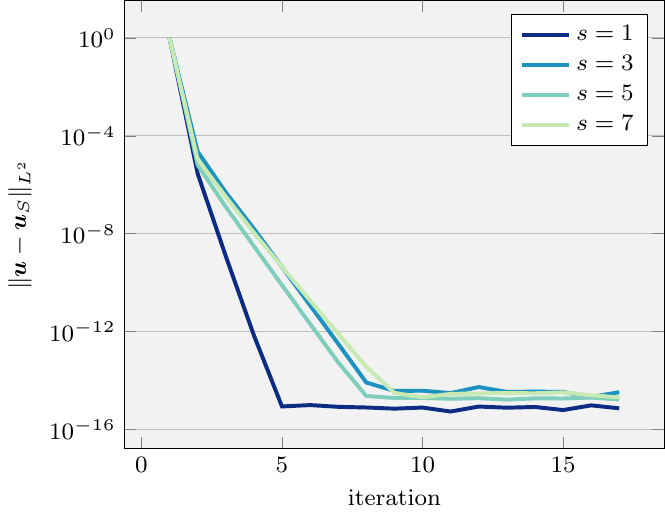}
  \end{minipage}                      
  \hfill
  \begin{minipage}[b]{0.49\textwidth}
    \centering
    \includegraphics{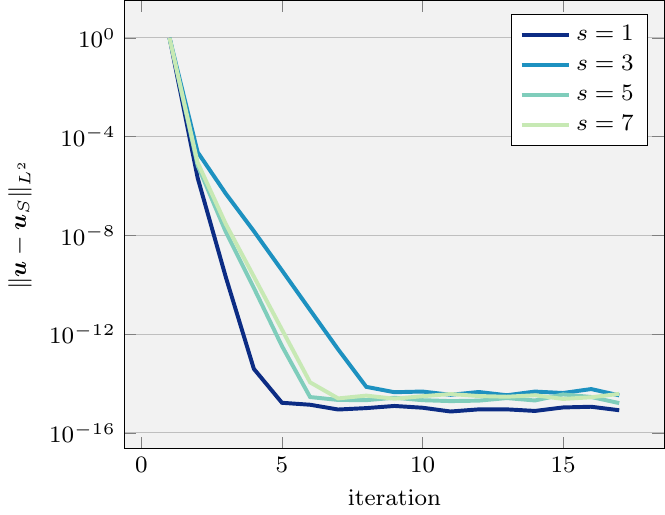}
  \end{minipage}
  \par
  \begin{minipage}[b]{0.49\textwidth}
    \centering
    \includegraphics{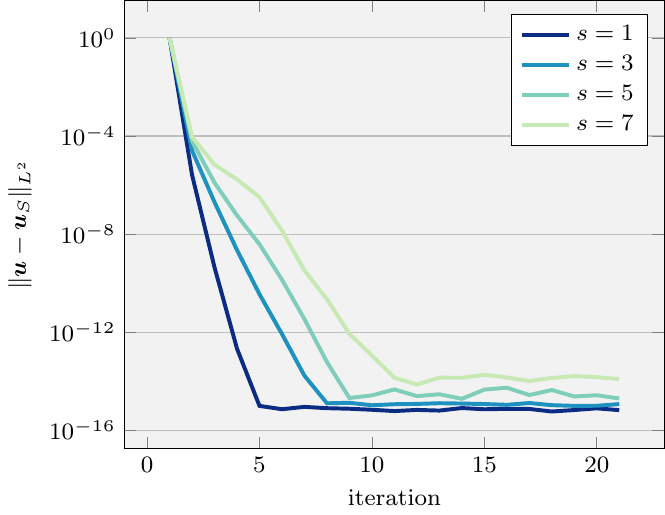}
  \end{minipage}                      
  \hfill
  \begin{minipage}[b]{0.49\textwidth}
    \centering
    \includegraphics{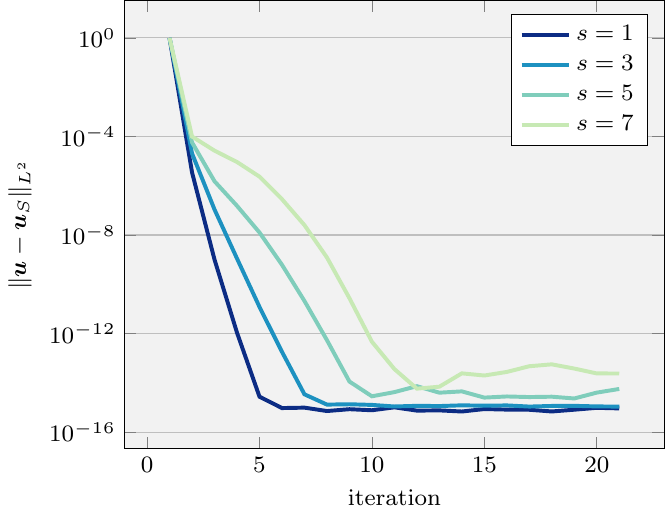}
  \end{minipage}
  \par
  \begin{minipage}[b]{0.49\textwidth}
    \centering
    \includegraphics{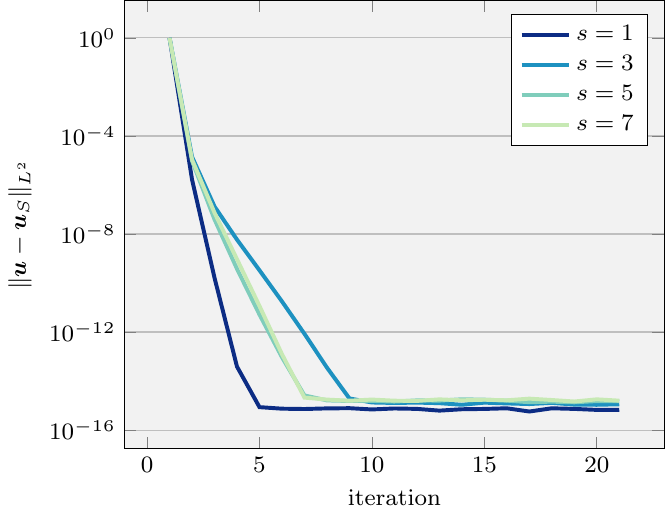}
  \end{minipage}                      
  \hfill
  \begin{minipage}[b]{0.49\textwidth}
    \centering
    \includegraphics{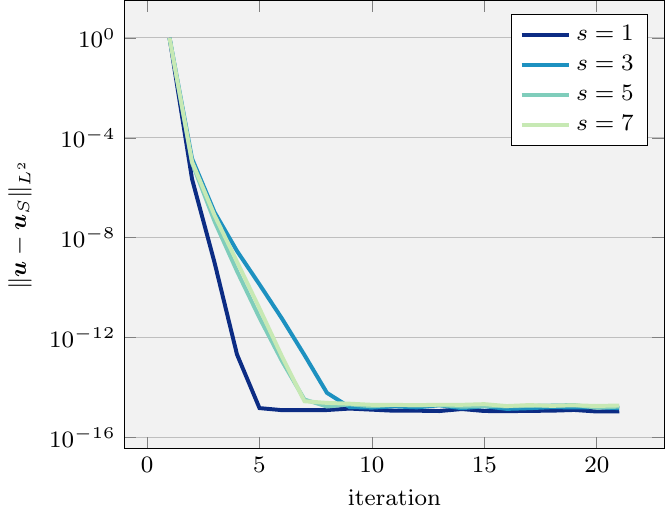}
  \end{minipage}
\fi
\caption[Convergence for high-order schemes]{Evolution of the error in the MGRIT iterates, applied to Burgers' (top), shallow-water (middle) and Euler equations (bottom). MGRIT parameters: V-cycle, F-smoothing, $3$ levels with $N_t=100,200,400$, respectively, $N_x=64$. Domain: $L=1$, $T=0.5$. Flux discretisation used: Lax-Friedrichs (left plots, \eqref{eqn::LFflux}), and Roe (right plots, \eqref{eqn::roeflux}), both employing WENO reconstruction of different orders of accuracy $s$. Time-stepper used: SSPRK3 \eqref{eqn::SSPRK3} with rediscretisation at all levels. Data generated running \texttt{testMGRIT\_highOrder.m} \cite{myCode}}
\label{fig::MGRITweno}
\end{figure*}
The initial conditions for the various problems are picked to ensure that the solutions achieve a similar maximum CFL number $c\sim 0.43$. These are:
\begin{equation}
  \begin{split}
    \bar{u}_B^0 =&\frac{4}{3}\sin\left(\frac{2\pi x}{L}\right),\; \bar{\boldsymbol{u}}_S^0=\frac{1}{11}\left[
    \begin{array}{cc}
      1+\frac{1}{2}\sin\left(\frac{2\pi x}{L}\right), & 0\\
    \end{array}\right]^T, \\
    \bar{\boldsymbol{u}}_E^0 =&\left[
    \begin{array}{ccc}
      1, & 0, & 1+\frac{1}{2}\sin\left(\frac{2\pi x}{L}\right)
    \end{array}\right]^T,
  \end{split}
\end{equation}
for Burgers', shallow-water, and Euler equations, respectively.
From \figref{fig::MGRITweno}, we see that the error behaviour varies in quite an erratic way as we modify $s$. For both Lax-Friedrichs and Roe fluxes with the SSPRK3 time-stepper, the best performance is seen with an integrator which completely disregards reconstruction, \emph{i.e.}, $s=1$. This is in line with expectations, as this choice of $s$ corresponds to the most dissipative among the schemes considered. The tendency is for performance to worsen as we increase $s$, although $s=3$ seems to provide an exception and behave particularly badly in some cases: the reason behind this remains unclear. Nonetheless, these higher-order methods still behave reasonably well, particularly in comparison to many results in the literature regarding more diffusive schemes. Evidence of this is also given in \figref{fig::MGRITcfl}: at least for $c\leq.95$, which corresponds to the dashed lines in the plots, SSPRK3 (top-right) outperforms FE (top-left), even if the latter is used in conjunction with the first order-matching procedure \eqref{eqn::coarseLFord1}. Indeed it is accuracy in the \emph{time} discretisation that helps convergence, while increasing the order of the \emph{space} discretisation seems to produce the opposite effect, as already discussed before, and further testified by the poorer error behaviour in the bottom plots of \figref{fig::MGRITcfl}. In particular, note how MGRIT (with rediscretisation on the coarse grid) fails to produce a converging solver for $c=.95$ (or even $c=0.48$) in the bottom-left of \figref{fig::MGRITcfl}, where it is used in combination with a high-order spatial discretisation and a low-order temporal discretisation.
\begin{figure*}[ht!]
\centering
\ifcompileFigs
  \begin{minipage}[b]{0.49\textwidth}
    \centering
    \begin{tikzpicture}[baseline,scale=0.8]
      \begin{semilogyaxis}[tick label style={font=\small},
                   xlabel= iteration,
                   ylabel= $\|\boldsymbol{u}-\boldsymbol{u}_S\|_{L^2}$,
                   ymajorgrids= true,
                   axis background/.style={fill=gray!10},
                   legend pos=south west,
                   cycle list = {{YlGnBu-F},{YlGnBu-L}} ]  
        \addplot+[line width=0.50mm,      ] table[x index=0, y index=5] {data/MGRIT_Burgers1DWENO_LF.dat};
        \addplot+[line width=0.50mm,      ] table[x index=0, y index=6] {data/MGRIT_Burgers1DWENO_LF.dat};
        \addplot+[line width=0.50mm,dashed] table[x index=0, y index=7] {data/MGRIT_Burgers1DWENO_LF.dat};
        \addplot+[line width=0.50mm,dashed] table[x index=0, y index=8] {data/MGRIT_Burgers1DWENO_LF.dat};
        \addlegendentry{$s\uparrow$ as $l\uparrow$}
        \addlegendentry{$s\downarrow$ as $l\uparrow$}
        \addlegendentry{$s,d\uparrow$ as $l\uparrow$}
        \addlegendentry{$s,d\downarrow$ as $l\uparrow$}
      \end{semilogyaxis}
    \end{tikzpicture}
  \end{minipage}                      
  \hfill
  \begin{minipage}[b]{0.49\textwidth}
    \centering
    \begin{tikzpicture}[baseline,scale=0.8]
      \begin{semilogyaxis}[tick label style={font=\small},
                   xlabel= iteration,
                   ylabel= $\|\boldsymbol{u}-\boldsymbol{u}_S\|_{L^2}$,
                   ymajorgrids= true,
                   axis background/.style={fill=gray!10},
                   legend pos=south west,
                   cycle list = {{YlGnBu-F},{YlGnBu-L}} ]  
        \addplot+[line width=0.50mm,      ] table[x index=0, y index=5] {data/MGRIT_Burgers1DWENO_ROE.dat};
        \addplot+[line width=0.50mm,      ] table[x index=0, y index=6] {data/MGRIT_Burgers1DWENO_ROE.dat};
        \addplot+[line width=0.50mm,dashed] table[x index=0, y index=7] {data/MGRIT_Burgers1DWENO_ROE.dat};
        \addplot+[line width=0.50mm,dashed] table[x index=0, y index=8] {data/MGRIT_Burgers1DWENO_ROE.dat};
        \addlegendentry{$s\uparrow$ as $l\uparrow$}
        \addlegendentry{$s\downarrow$ as $l\uparrow$}
        \addlegendentry{$s,d\uparrow$ as $l\uparrow$}
        \addlegendentry{$s,d\downarrow$ as $l\uparrow$}
      \end{semilogyaxis}
    \end{tikzpicture}
  \end{minipage}
\else
\begin{minipage}[b]{0.49\textwidth}
    \centering
    \includegraphics{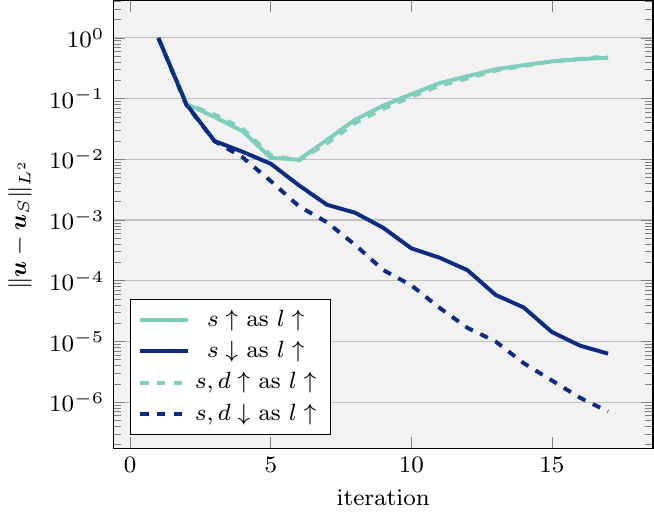}
  \end{minipage}                      
  \hfill
  \begin{minipage}[b]{0.49\textwidth}
    \centering
    \includegraphics{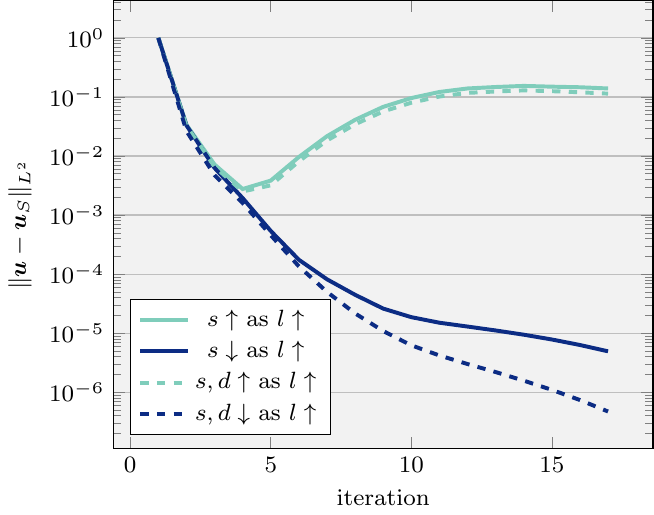}
  \end{minipage}
\fi
\caption[Convergence when varying accuracy across levels]{Impact of modifying the accuracy of the schemes used across the levels on the error convergence of MGRIT.
Same parameters as in \figref{fig::MGRITweno}, but with the order of the WENO reconstruction $s$ either increasing (green) or decreasing (blue) as we move to coarser levels. Dashed lines show error convergence if both $s$ and the order accuracy of the time-stepper, $d$, are changed simultaneously, while solid lines show results if only $s$ is varied. Left is for Lax-Friedrichs flux, right is for Roe flux. Only results for the application to Burgers' equation are shown, but the ones collected for the other test-cases follow the same qualitative behaviour. Data generated via \texttt{testMGRIT\_highOrder.m} \cite{myCode}}
\label{fig::MGRITweno_diffOrd}
\end{figure*}

Ultimately, however, even in the high-order case the CFL condition has shown to be the true bottleneck in the effectiveness of such solvers. Still looking at the plots in \figref{fig::MGRITcfl}, we can see that high-order schemes are much more sensitive to an increase in $c$. This clearly highlights that, if explicit solvers are to be used, limiting ourselves to coarsening just in time, and to employing coarse integrators based on rediscretisation, is not a viable option: an effective algorithm needs to pair coarsening in space and time together, in order to keep the CFL number close to $1$ at each level.

\paragraph{Changing accuracy}
Experimenting with discretisations of different accuracy at different multigrid levels makes sense in attempting to match the action of the fine solver: one can think of making up for the loss of precision due to coarsening by increasing the order of the reconstruction at a coarser level, and still obtain a result similar to what the fine integrator would provide. With this motivation in mind, we test applications of MGRIT where the order of the discretisations used, (both in space and time), varies across the levels.
Unfortunately, the results reported in \figref{fig::MGRITweno_diffOrd} pinpoint that this strategy is not beneficial to the algorithm performance and that, instead, simple rediscretisation is preferable. Among the alternatives considered, though, decreasing the accuracy as we descend to coarser levels has proven the most promising, which could make it attractive if one seeks to make the application of coarse integrators cheaper: this is also in line with the results shown in \cite{Nielsen}, where a lower-order WENO reconstruction was used at the coarse level.

\section{Conclusion}
In this paper, we consider the use of MGRIT for the parallelisation in time of the solution of non-linear hyperbolic PDEs. In particular, we aim to understand how the choice of integrators used at each level of the multigrid algorithm impacts its convergence behaviour. To this purpose, we have measured the performance of MGRIT applied to a number of test problems, using a combination of existing integrators commonly used for the solution of conservation laws. The results show that dissipative schemes behave better in general, which is in line with the literature; however, higher-order methods for spatial reconstruction coupled with matching order time-discretisations still provide satisfactory convergence results, under suitable CFL limits.

We note the importance of choosing coarse integrators that closely follow the action of the fine integrator, by showcasing the effectiveness of a new method proposed for the construction of accurate coarse solvers. This approach seeks to directly approximate the action of the fine solver: even though it comes at an increased cost with respect to simple rediscretisation, it offers superior performance. Its range of applicability remains, however, limited.

In all cases, the performance is seen to degrade fairly quickly as the CFL number increases. This seems to be a strong limitation for the application of MGRIT to hyperbolic systems directly discretised with explicit time-steppers. The level of coarsening that can be applied to the temporal grid in this case is, thus, effectively capped. This suggests that MGRIT should be paired with a spatial coarsening strategy as well, in order to control the CFL number adequately across all levels.

\bibliographystyle{siamplain}      
\bibliography{main}   				

\begin{thebibliography}{10}

\bibitem{pararealAnalysis}
{\sc G.~Bal}, {\em On the convergence and the stability of the parareal
  algorithm to solve partial differential equations}, in Domain Decomposition
  Methods in Science and Engineering, T.~J. Barth, M.~Griebel, D.~E. Keyes,
  R.~M. Nieminen, D.~Roose, T.~Schlick, R.~Kornhuber, R.~Hoppe, J.~P{\'e}riaux,
  O.~Pironneau, O.~Widlund, and J.~Xu, eds., Berlin, Heidelberg, 2005, Springer
  Berlin Heidelberg, pp.~425--432.

\bibitem{BoltenEtAl2018}
{\sc M.~Bolten, D.~Moser, and R.~Speck}, {\em Asymptotic convergence of the
  parallel full approximation scheme in space and time for linear problems},
  Numerical Linear Algebra with Applications, 25 (2018), p.~e2208,
  \url{https://doi.org/10.1002/nla.2208}.

\bibitem{multigrid}
{\sc W.~Briggs, V.~Henson, and S.~McCormick}, {\em A Multigrid Tutorial: Second
  Edition}, Other Titles in Applied Mathematics, Society for Industrial and
  Applied Mathematics, 2000.

\bibitem{RIDC}
{\sc A.~Christlieb, C.~B.~MacDonald, B.~W.~Ong, and R.~Spiteri}, {\em
  Revisionist integral deferred correction with adaptive stepsize control},
  Communications in Applied Mathematics and Computational Science, 10 (2013),
  \url{https://doi.org/10.2140/camcos.2015.10.1}.

\bibitem{ChristliebEtAl2010}
{\sc A.~J. Christlieb, C.~B. Macdonald, and B.~W. Ong}, {\em {Parallel
  high-order integrators}}, SIAM Journal on Scientific Computing, 32 (2010),
  pp.~818--835, \url{https://doi.org/10.1137/09075740X}.

\bibitem{pararealHyp}
{\sc X.~Dai and Y.~Maday}, {\em Stable parareal in time method for first- and
  second-order hyperbolic systems}, SIAM Journal on Scientific Computing, 35
  (2013), pp.~A52--A78, \url{https://doi.org/10.1137/110861002}.

\bibitem{myCode}
{\sc F.~Danieli}, {\em Repo {MGRIT\_N}on-linear{\_}hyperbolic}.
\newblock \url{https://gitlab.com/fdanieli/mgrit_non-linear_hyperbolic}.

\bibitem{MGRITconv}
{\sc V.~Dobrev, T.~Kolev, N.~Petersson, and J.~Schroder}, {\em Two-level
  convergence theory for multigrid reduction in time {(MGRIT)}}, SIAM Journal
  on Scientific Computing, 39 (2017), pp.~S501--S527,
  \url{https://doi.org/10.1137/16M1074096}.

\bibitem{pararealHyp2}
{\sc A.~Eghbal, A.~G. Gerber, and E.~Aubanel}, {\em Acceleration of unsteady
  hydrodynamic simulations using the parareal algorithm}, Journal of
  Computational Science, 19 (2017), pp.~57 -- 76,
  \url{https://doi.org/https://doi.org/10.1016/j.jocs.2016.12.006}.

\bibitem{PFASST}
{\sc M.~Emmett and M.~L. Minion}, {\em {Toward an Efficient Parallel in Time
  Method for Partial Differential Equations}}, Communications in Applied
  Mathematics and Computational Science, 7 (2012), pp.~105--132,
  \url{https://doi.org/10.2140/camcos.2012.7.105}.

\bibitem{MGRIToriginal}
{\sc R.~Falgout, S.~Friedhoff, T.~Kolev, S.~MacLachlan, and J.~Schroder}, {\em
  Parallel time integration with multigrid}, SIAM Journal on Scientific
  Computing, 36 (2014), pp.~C635--C661,
  \url{https://doi.org/10.1137/130944230}.

\bibitem{MGRITNonLin}
{\sc R.~Falgout, T.~Manteuffel, B.~O'Neill, and J.~Schroder}, {\em Multigrid
  reduction in time for nonlinear parabolic problems: A case study}, SIAM
  Journal on Scientific Computing, 39 (2017), pp.~S298--S322,
  \url{https://doi.org/10.1137/16M1082330}.

\bibitem{MGRIT_BDF}
{\sc R.~D. Falgout, S.~Friedhoff, T.~V. Kolev, S.~P. MacLachlan, J.~B.
  Schroder, and S.~Vandewalle}, {\em Multigrid methods with space--time
  concurrency}, Computing and Visualization in Science, 18 (2017),
  pp.~123--143, \url{https://doi.org/10.1007/s00791-017-0283-9}.

\bibitem{MGRIT_BDF2}
{\sc R.~D. Falgout, M.~Lecouvez, and C.~S. Woodward}, {\em A parallel-in-time
  algorithm for variable step multistep methods}, in LLNL Technical Report.

\bibitem{pararealAnalysis2}
{\sc M.~Gander and S.~Vandewalle}, {\em Analysis of the parareal time-parallel
  time-integration method}, SIAM Journal on Scientific Computing, 29 (2007),
  pp.~556--578, \url{https://doi.org/10.1137/05064607X}.

\bibitem{50yr}
{\sc M.~J. Gander}, {\em 50 years of time parallel time integration}, in
  Multiple Shooting and Time Domain Decomposition Methods, T.~Carraro,
  M.~Geiger, S.~K\"orkel, and R.~Rannacher, eds., Springer International
  Publishing, 2014, ch.~3, pp.~69--113.

\bibitem{SSPRK}
{\sc S.~Gottlieb, C.~Shu, and E.~Tadmor}, {\em Strong stability-preserving
  high-order time discretization methods}, SIAM Review, 43 (2001), pp.~89--112,
  \url{https://doi.org/10.1137/S003614450036757X}.

\bibitem{HHentropyFix}
{\sc A.~Harten and J.~M. Hyman}, {\em Self adjusting grid methods for
  one-dimensional hyperbolic conservation laws}, Journal of Computational
  Physics, 50 (1983), pp.~235 -- 269,
  \url{https://doi.org/https://doi.org/10.1016/0021-9991(83)90066-9}.

\bibitem{Beth}
{\sc T.~Haut and B.~Wingate}, {\em An asymptotic parallel-in-time method for
  highly oscillatory {PDEs}}, SIAM Journal on Scientific Computing, 36 (2014),
  pp.~A693--A713, \url{https://doi.org/10.1137/130914577}.

\bibitem{jan}
{\sc J.~Hesthaven}, {\em Numerical Methods for Conservation Laws}, Society for
  Industrial and Applied Mathematics, Philadelphia, PA, 2017,
  \url{https://doi.org/10.1137/1.9781611975109}.

\bibitem{pararealSkinTransport}
{\sc A.~Kreienbuehl, A.~Naegel, D.~Ruprecht, R.~Speck, G.~Wittum, and
  R.~Krause}, {\em Numerical simulation of skin transport using parareal},
  Computing and Visualization in Science, 17 (2015), pp.~99--108,
  \url{https://doi.org/10.1007/s00791-015-0246-y}.

\bibitem{leveque}
{\sc R.~LeVeque}, {\em Numerical Methods for Conservation Laws}, Lectures in
  Mathematics ETH Z{\"u}rich, Department of Mathematics Research Institute of
  Mathematics, Springer, 1992.

\bibitem{parareal}
{\sc J.-L. Lions, Y.~Maday, and G.~Turinici}, {\em {R{\'e}solution d'EDP par un
  sch{\'e}ma en temps {\tt<\tt<}parar{\'e}el{\tt>\tt>}}}, {Comptes rendus de
  l'Acad{\'e}mie des sciences. S{\'e}rie I, Math{\'e}matique}, 332 (2001),
  pp.~661--668.

\bibitem{WENO}
{\sc X.-D. Liu, S.~Osher, and T.~Chan}, {\em Weighted essentially
  non-oscillatory schemes}, Journal of Computational Physics, 115 (1994),
  pp.~200 -- 212, \url{https://doi.org/https://doi.org/10.1006/jcph.1994.1187}.

\bibitem{Nielsen}
{\sc A.~S. Nielsen, G.~Brunner, and J.~S. Hesthaven}, {\em Communication-aware
  adaptive parareal with application to a nonlinear hyperbolic system of
  partial dierential equations}, Journal of Computational Physics, 15 (2017),
  pp.~483--505, \url{https://doi.org/10.1016/j.jcp.2018.04.056}.

\bibitem{OngEtAl2016}
{\sc B.~W. Ong, R.~D. Haynes, and K.~Ladd}, {\em Algorithm 965: {RIDC} methods:
  A family of parallel time integrators}, ACM Trans. Math. Softw., 43 (2016),
  pp.~8:1--8:13, \url{https://doi.org/10.1145/2964377}.

\bibitem{pararealOptionPricing}
{\sc G.~Pagès, O.~Pironneau, and G.~Sall}, {\em The parareal algorithm for
  {A}merican options}, Comptes Rendus Mathematique, 354 (2016), pp.~1132 --
  1138, \url{https://doi.org/https://doi.org/10.1016/j.crma.2016.09.010}.

\bibitem{WENOcharRec}
{\sc J.~Qiu and C.-W. Shu}, {\em On the construction, comparison, and local
  characteristic decomposition for high-order central {WENO} schemes}, Journal
  of Computational Physics, 183 (2002), pp.~187 -- 209,
  \url{https://doi.org/https://doi.org/10.1006/jcph.2002.7191}.

\bibitem{processorPower}
{\sc K.~Rupp}, {\em 42 years of microprocessor trend data}.
\newblock
  \url{https://www.karlrupp.net/2018/02/42-years-of-microprocessor-trend-data/}.
\newblock Raw data available at
  \url{https://github.com/karlrupp/microprocessor-trend-data}. Accessed on:
  2019/06/19.

\bibitem{pararealAnalysis3}
{\sc D.~Ruprecht}, {\em Wave propagation characteristics of parareal},
  Computing and Visualization in Science, 19 (2018), pp.~1--17,
  \url{https://doi.org/10.1007/s00791-018-0296-z}.

\bibitem{WENO2}
{\sc C.-W. Shu}, {\em Essentially non-oscillatory and weighted essentially
  non-oscillatory schemes for hyperbolic conservation laws}, Springer Berlin
  Heidelberg, Berlin, Heidelberg, 1998, pp.~325--432,
  \url{https://doi.org/10.1007/BFb0096355}.

\bibitem{MGRIT_Oliver}
{\sc H.~D. Sterck, R.~D. Falgout, S.~Friedhoff, O.~A. Krzysik, and S.~P.
  MacLachlan}, {\em Pseudo-optimal parareal and {MGRIT} coarse-grid operators
  for linear advection}.
\newblock in preparation.

\bibitem{ToselliWidlundDDM}
{\sc A.~Toselli and O.~Widlund}, {\em Domain Decomposition Methods - Algorithms
  and Theory}, Springer Series in Computational Mathematics, Springer Berlin
  Heidelberg, 2006.

\end{thebibliography}


\end{document}